\documentclass[12pt,leqno]{article}
\usepackage{amssymb}
\usepackage[mathscr]{eucal}
\usepackage{amsmath,amssymb,latexsym,theorem,enumerate,ifthen}
\usepackage[english]{babel}
\usepackage{color,url}
\usepackage{appendix}

\newcommand{\NN}{\mathbb{N}}
\newcommand{\QQ}{\mathbb{Q}}
\newcommand{\RR}{\mathbb{R}}

\newcommand{\ZZ}{\mathbb{Z}}

\newcommand{\bone}{{\boldsymbol{1}}}

\newcommand{\cA}{{\mathcal A}}
\newcommand{\cB}{{\mathcal B}}

\newcommand{\cD}{{\mathcal D}}

\newcommand{\cG}{{\mathcal G}}
\newcommand{\cF}{{\mathcal F}}

\newcommand{\cN}{{\mathcal N}}

\newcommand{\cX}{{\mathcal X}}

\newcommand{\dd}{\mathrm{d}}
\newcommand{\ee}{\mathrm{e}}

\newcommand{\EE}{\operatorname{\mathbb{E}}}

\newcommand{\PP}{\operatorname{\mathbb{P}}}

\newcommand{\var}{\operatorname{Var}}
\newcommand{\cov}{\operatorname{Cov}}
\newcommand{\sign}{\operatorname{sign}}

\newcommand{\esssup}{\operatorname{ess\,sup}}
\DeclareMathOperator*{\argmin}{arg\,min}

\newcommand{\vare}{\varepsilon}

\renewcommand{\mid}{\,|\,}

\renewcommand{\leq}{\leqslant}
\renewcommand{\geq}{\geqslant}

\newcommand{\stoch}{\stackrel{\PP}{\longrightarrow}}
\newcommand{\distr}{\stackrel{\cD}{\longrightarrow}}
\newcommand{\distre}{\stackrel{\cD}{=}}

\newcommand{\as}{\stackrel{{\mathrm{a.s.}}}{\longrightarrow}}
\newcommand{\ase}{\stackrel{{\mathrm{a.s.}}}{=}}

\newcommand{\proofend}{\hfill\mbox{$\Box$}}

\numberwithin{equation}{section}

\theoremstyle{change} \theorembodyfont{\em}
\newtheorem{Lem}{Lemma.}[section]
\newtheorem{Thm}[Lem]{Theorem.}

\newtheorem{Cor}[Lem]{Corollary.}
\newtheorem{Def}[Lem]{Definition.}

\newtheorem{claim}[Lem]{Claim.}

\theorembodyfont{\rm}
\newtheorem{Rem}[Lem]{Remark.}
\newtheorem{Ex}[Lem]{Example.}

\def\OnlyOnArXiv#1#2{\ifthenelse{\equal{#1}{Y}}{#2}{}}

\begin{document}

\begin{center}
 {\bfseries\Large Limit theorems for deviation means of independent and identically distributed random variables}

\vspace*{3mm}

{\sc\large
  M\'aty\'as $\text{Barczy}^{*,\diamond}$,
  Zsolt $\text{P\'ales}^{**}$ }

\end{center}

\vskip0.2cm

\noindent
 * MTA-SZTE Analysis and Stochastics Research Group,
   Bolyai Institute, University of Szeged,
   Aradi v\'ertan\'uk tere 1, H--6720 Szeged, Hungary.

\noindent
 ** Institute of Mathematics, University of Debrecen,
    Pf.~400, H--4002 Debrecen, Hungary.

\noindent e-mail: barczy@math.u-szeged.hu (M. Barczy),
                  pales@science.unideb.hu  (Zs. P\'ales).

\noindent $\diamond$ Corresponding author.

\vskip0.2cm



{\renewcommand{\thefootnote}{}
\footnote{\textit{2020 Mathematics Subject Classifications\/}:
 60F05, 60F15, 26E60, 62F12}
\footnote{\textit{Key words and phrases\/}:
 deviation mean, Bajraktarevi\'c mean, strong law of large numbers, central limit theorem, law of the iterated logarithm, large deviations.}
\vspace*{0.2cm}
\footnote{M\'aty\'as Barczy was supported by the project TKP2021-NVA-09.
Project no.\ TKP2021-NVA-09 has been implemented with the support
 provided by the Ministry of Innovation and Technology of Hungary from the National Research, Development and Innovation Fund,
 financed under the TKP2021-NVA funding scheme.
Zsolt P\'ales is supported by the K-134191 NKFIH Grant.}}

\vspace*{-10mm}

\begin{abstract}
We derive a strong law of large numbers, a central limit theorem, a law of the iterated logarithm and a large deviation theorem for so-called deviation means of independent and identically distributed random variables (for the strong law of large numbers, we suppose only pairwise independence instead of (total) independence).
The class of deviation means is a special class of M-estimators or more generally extremum estimators,
  which are well-studied in statistics.
The assumptions of our limit theorems for deviation means seem to be new and weaker than the known ones for M-estimators in the literature.
Especially, our results on the strong law of large numbers and on the central limit theorem generalize
 the corresponding ones for quasi-arithmetic means due to de Carvalho (2016) and the ones for Bajraktarevi\'c means
 due to Barczy and Burai (2022).
\end{abstract}

\section{Introduction}
\label{section_intro}

The theory and applications of means can be traced back to the ancient Greeks who developed the notions of arithmetic, geometric and harmonic means.
Since then various kinds of means have been introduced and their analytical and algebraic properties have been systemically studied by many researchers,
 for recent surveys, see e.g., Bullen \cite{Bul03} and Beliakov et al.\ \cite{BelBusCal16}.
Studying probabilistic properties of means got less attention in general.
In the present paper we investigate such a property of deviation means (see Definition \ref{Def_deviation_mean}) introduced by Dar\'oczy \cite{Dar71b,Dar72b},
  namely, we study whether deviation means have a centralizing property such as the usual arithmetic mean
  in the strong law of large numbers, the central limit theorem, the law of the iterated logarithm and the theory of large deviations for independent and
  identically distributed (i.i.d.) random variables.
The class of deviation means is in fact a special class of M-estimators introduced by Huber \cite{Hub64, Hub67},
 where the letter M refers to "maximum likelihood-type".
More precisely, deviation means are special cases of so-called $\psi$-estimators (also known as Z-estimators).
In Remark \ref{Rem_M_estimate} we recall M-estimators and then in Remark \ref{Rem_M_estimate_deviation_mean} we present deviation means as special M-estimators.
Here we only note that it seems to us that the theory of deviation means and that of M-estimators developed independently of each other.
In Remark \ref{Rem_M_estimate} we also point out that M-estimators are in fact special extremum estimators (see, e.g., Hayashi \cite[Section 7.1]{Hay}).

In statistics there is a huge literature on limit theorems for M-estimators and extremum estimators
 based on i.i.d.\ random variables and on more general observations.
Here our goal is not to present the development of M-estimators and extremum estimators or to give a full bibliography of contributions to the literature.
Rather, we want to focus on limit theorems in case of the special class of deviation means, and to prove limit theorems
 under conditions which seem to be new and weaker then the known ones for M-estimators in the literature.
In Remark \ref{Rem_feltet_osszehas} we give a comparison of the conditions for our limit theorems with those of some of the existing results for M-estimators.
There are plenty sets of sufficient assumptions under which limit theorems for M-estimators hold,
 and in Remark \ref{Rem_feltet_osszehas} we also point out that our set of sufficient conditions for deviation means
 are different from the ones available in the literature, and we think that they are more simpler and more checkable.
For quasi-arithmetic means (see Definition \ref{Def_quasi_arithmetic}) and Bajraktarevi\'c means (see Definition \ref{Def_Bajraktarevic}), which are special deviation means, strong laws of large numbers and central limit theorems have already been derived by de Carvalho \cite{Car} and Barczy and Burai \cite{BarBur} without referring to general limit theorems for M-estimators and extremum estimators.

First, we recall some related literature on limit theorems for some other kinds of means of i.i.d.\ random variables.
The asymptotic behaviour of a special mean related to the elementary symmetric polynomials of i.i.d.\ positive random variables attracted several researchers' attention.
Namely, for positive integers $k$, $n$ with $k\leq n$ and positive real numbers $x_1,\ldots,x_n$, let
 \[
   M_n^{(k)}(x_1,\ldots,x_n):=\sqrt[k]{ \frac{\sum_{1\leq i_1<i_2<\cdots< i_k\leq n}x_{i_1}\cdots x_{i_k}}{\binom{n}{k}}},
 \]
 which is an $n$-variable mean in the sense of the forthcoming Definiton \ref{Def_mean}.
This mean could be called Newton's mean, since Newton proved that for positive integers $k,n$ with $2\leq k\leq n-1$ and positive real numbers
 $x_1,\ldots,x_n$, we have
 \[
   M_n^{(k-1)}(x_1,\ldots,x_n) M_n^{(k+1)}(x_1,\ldots,x_n) \leq (M_n^{(k)}(x_1,\ldots,x_n))^2,
 \]
 see, e.g., Hardy et al.\ \cite[Theorems 51 and 144]{HarLitPol34}. This inequality easily yields that
 \[
  M_n^{(1)}(x_1,\ldots,x_n)\geq M_n^{(2)}(x_1,\ldots,x_n) \geq\cdots\geq M_n^{(n)}(x_1,\ldots,x_n),
   \qquad x_1,\ldots,x_n>0,
 \]
 where, one can observe that $M_n^{(1)}(x_1,\ldots,x_n)$ and $M_n^{(n)}(x_1,\ldots,x_n)$ are nothing else but the arithmetic and geometric means of $x_1,\ldots,x_n$, respectively, see Hardy et al.\ \cite[Theorem 52]{HarLitPol34}.
Let us given a sequence $(\xi_n)_{n\geq 1}$ of i.i.d.\ positive random variables.
Under some appropriate moment conditions on $\xi_1$, Hoeffding \cite{Hoe} and Major \cite[Theorem 1]{Maj} proved central limit theorems
 for $(M_n^{(k)}(\xi_1,\ldots,\xi_n))^k$ and $\ln((M_n^{(k)}(\xi_1,\ldots,\xi_n))^k)+\ln(\binom{n}{k})$, respectively, as $n\to\infty$;
 Hal\'asz and Sz\'ekely \cite{HalSzek} established a strong law of large numbers for $M_n^{(k)}(\xi_1,\ldots,\xi_n)$ as $n\to\infty$
 (in this case $k$ may depend on $n$ as well); and Rempala and Gupta \cite[Theorems 2 and 4]{RemGup} derived laws of the iterated logarithm
  for  $(M_n^{(k)}(\xi_1,\ldots,\xi_n))^k$ and  $\ln((M_n^{(k)}(\xi_1,\ldots,\xi_n))^k)$ as $n\to\infty$.

Let $\NN$, $\ZZ_+$, $\QQ$, $\RR$ and $\RR_+$ denote the sets of positive integers, non-negative integers, rational numbers,
real numbers and non-negative real numbers. An interval $I\subset\RR$ will be called nondegenerate if it contains at least two distinct points. All the random variables are defined on a probability space $(\Omega,\cF,\PP)$.
Convergence almost surely, in probability and in distribution will be denoted by $\as$, $\stoch$ and $\distr$, respectively. The almost surely equality and the equlity in distribution is denoted by $\ase$ and $\distre$, respectively.
The cumulative distribution function of a random variable $\eta:\Omega\to\RR$ is defined by $F_\eta:\RR\to [0,1]$, $F_\eta(y):=\PP(\eta<y)$, $y\in\RR$.
For a random variable $\eta$, its essential supremum is defined by
 \[
   \esssup(\eta) := \sup\Big\{x\in\RR : \PP(\eta\leq x)<1\Big\}.
 \]
For any $\sigma>0$, $\cN(0,\sigma^2)$ denotes a normal distribution with mean $0$ and variance $\sigma^2$.

\begin{Def}\label{Def_mean}
Let $I$ be a nondegenerate interval of $\RR$ and $n\in\NN$.
A function $M:I^n \to \RR$ is called an $n$-variable mean on $I$ if
 \begin{align}\label{def_mean}
    \min(x_1,\ldots,x_n) \leq M(x_1,\ldots,x_n) \leq \max(x_1,\ldots,x_n),
     \qquad x_1,\ldots,x_n\in I.
 \end{align}
 An $n$-variable mean $M$ on $I$ is called strict if both inequalities in \eqref{def_mean} are sharp for all $x_1,\ldots,x_n\in I$ \
 satisfying $\min(x_1,\ldots,x_n)< \max(x_1,\ldots,x_n)$.
\end{Def}

If $n=1$, then the only $1$-variable mean $M$ on $I$ is $M(x)=x$, $x\in I$.

A classical and well-studied class of means is the class of quasi-arithmetic means (see, e.g., the monograph of Hardy et al.\ \cite{HarLitPol34}).

\begin{Def}[Quasi-arithmetic mean]\label{Def_quasi_arithmetic}
Let $n\in\NN$, let $I$ be a nondegenerate interval  of $\RR$, and let $f:I\to\RR$ be a continuous and strictly increasing function.
The $n$-variable quasi-arithmetic mean $\mathscr{A}^f_n:I^n\to I$ is defined by
 \[
  \mathscr{A}^f_n(x_1,\ldots,x_n):= f^{-1} \bigg( \frac{1}{n} \sum_{i=1}^n f(x_i)\bigg), \qquad x_1,\ldots, x_n\in I,
 \]
where $f^{-1}$ denotes the inverse of $f$. The function $f$ is called the generator of $\mathscr{A}^f_n$.
\end{Def}

Kolmogorov \cite{Kol30}, Nagumo \cite{Nag30, Nag31}, and de Finetti \cite{Def31}
 (see also Tikhomirov \cite[page 144]{Tik}), independently of each other, established the following axiomatic characterization of quasi-arithmetic means.

\begin{Thm}[Kolmogorov (1930), Nagumo (1930) and de Finetti (1931)]\label{Thm_Kolmogorov}
Let $I$ be a compact nondegenerate interval of $\RR$ and let $M_n:I^n\to\RR$, $n\in\NN$, be a sequence of functions. Then the following two statements are equivalent:
\vspace{-8pt}
\begin{enumerate}[(i)]\itemsep=-4pt
\item The sequence $(M_n)_{n\in\NN}$ is quasi-arithmetic, that is, there exists a continuous and strictly monotone function $f:I\to\RR$ such that
 \[
   M_n(x_1,\ldots,x_n) = \mathscr{A}^f_n(x_1,\ldots,x_n),
    \qquad x_1,\ldots,x_n\in I,\;\; n\in\NN.
 \]
\item The sequence $(M_n)_{n\in\NN}$ possesses the following properties:
\begin{itemize}
  \item  $M_n$ is continuous and strictly increasing in each variable for each $n\in\NN$,
  \item $M_n$ is symmetric for each $n\in\NN$  (i.e., $M_n(x_1,\ldots,x_n) = M_n(x_{\pi(1)},\ldots,x_{\pi(n)})$ for each
        $x_1,\ldots,x_n \in I$ and each permutation $(\pi(1),\ldots,\pi(n))$ of $(1,\ldots,n)$),
  \item $M_n(x_1,\ldots,x_n) = x$ whenever $x_1=\cdots =x_n = x\in I$, $n\in\NN$,
  \item $M_{n+m}(x_1,\ldots,x_n, y_1,\ldots,y_m) = M_{n+m}(\overline x_n, \ldots, \overline x_n, y_1,\ldots, y_m)$
        for each $n,m\in\NN$, $x_1,\ldots,x_n,y_1,\ldots,y_m\in I$, where $\overline x_n:=M_n(x_1,\ldots,x_n)$.
 \end{itemize}
\end{enumerate}
\end{Thm}

\begin{Rem}
(i).
The arithmetic, geometric and harmonic mean are quasi-arithmetic means corresponding to the functions $f:\RR\to \RR$, $f(x):=x$, $x\in\RR$;
 $f:(0,\infty)\to \RR$, $f(x):=\ln(x)$, $x>0$; and $f:(0,\infty)\to \RR$, $f(x) =x^{-1}$, $x>0$, respectively.

(ii). Two quasi-arithmetic means on $I$, generated by $f$ and $g$, are equal (i.e., $\mathscr{A}^f_n=\mathscr{A}^g_n$ on $I^n$ for each $n\in\NN$)
 if and only if there exist $a,b\in\RR$, $a\ne 0$ such that
 \begin{align}\label{help_quasi_generator}
  f(x)=ag(x)+b, \qquad x\in I,
 \end{align}
see, e.g., Hardy et al.\ \cite[Section 3.2, Theorem 83]{HarLitPol34}.
As a consequence, the function $f$ in part (i) of Theorem \ref{Thm_Kolmogorov} can be chosen to be strictly increasing as well.
\proofend
\end{Rem}

Generalizing the notion of quasi-arithmetic means, Bajraktarevi\'c \cite{Baj58} introduced a new class of means
 (nowadays called Bajraktarevi\'c means) in the following way.

\begin{Def}[Bajraktarevi\'c mean]\label{Def_Bajraktarevic}
Let $n\in\NN$, let $I$ be a nondegenerate interval of $\RR$,
 let $f:I\to\RR$ be a continuous and strictly monotone
 function, and let $p:I\to (0,\infty)$ be a (weight) function.
The $n$-variable Bajraktarevi\'c mean $\mathscr{B}^{f,p}_n:I^n\to I$ is defined by
 \[
  \mathscr{B}^{f,p}_n(x_1,\ldots,x_n):= f^{-1}\left( \frac{\sum_{i=1}^n p(x_i) f(x_i)}{\sum_{i=1}^n p(x_i)} \right),\qquad x_1,\ldots, x_n\in I.
 \]
 The pair $(f,p)$ is called the generator of $\mathscr{B}^{f,p}_n$.
\end{Def}

For each $n\in\NN$, $\mathscr{B}^{f,p}_n$ is a strict, symmetric $n$-variable mean, see, e.g., Bajraktarevi\'c \cite{Baj58} or P\'ales and Zakaria \cite{PalZak20a}.
By choosing $p(x)=1$, $x\in I$, one can see that the Bajraktarevi\'c mean of $x_1,\ldots, x_n$ corresponding to $f$ and $p$ coincides
 with the quasi-arithmetic mean of $x_1,\ldots,x_n$ generated by $f$.
For the equality problem of Bajraktarevi\'c means, we refer to Losonczi \cite{Los99, Los06} and Losonczi et al.\ \cite[Theorem 16]{LosPalZak}.
The characterization theorem of Bajraktarevi\'c means was established by P\'ales \cite{Pal87}.

Next we recall the notions of a deviation and a deviation mean introduced by Dar\'oczy \cite{Dar71b,Dar72b}.

\begin{Def}[Deviation]\label{Def_deviation}
Let $I$ be a nondegenerate interval of $\RR$.
A function $D:I^2\to\RR$ is called a deviation if
 \vspace{-12pt}
 \begin{itemize}\itemsep=-4pt
   \item[(i)] for all $x\in I$, the function $I\ni t\mapsto D(x,t)$ is strictly decreasing and continuous,
   \item[(ii)] $D(t,t)=0$, $t\in I$.
 \end{itemize}
\end{Def}

\begin{Def}[Deviation mean]\label{Def_deviation_mean}
Let $n\in\NN$, let $I$ be a nondegenerate interval of $\RR$,  and let $D:I^2\to\RR$ be a deviation.
The $n$-variable deviation mean of $x_1,\ldots, x_n\in I$ corresponding to the deviation $D$ is defined as the unique solution of the equation
 \[
   \sum_{i=1}^n D(x_i,t) = 0, \qquad t\in I,
 \]
and is denoted by $\mathscr{M}_n^{D}(x_1,\ldots,x_n)$. The deviation $D$ is called the generator of $\mathscr{M}_n^{D}:I^n\to I$.
\end{Def}

Here we call the attention to the fact that in the expression 'deviation mean' the word 'deviation' has nothing to do with the word 'deviation' in the expression 'large deviations', which is an important subfield in probability theory.

In the next remark we explain that the notion of a deviation mean is well-defined.

\begin{Rem}\label{Rem1}
Let $I$ be a nondegenerate interval of $\RR$, and let $D:I^2\to\RR$ be a deviation.
Given $n\in\NN$ and $x_1,\dots,x_n\in I$, by (i) of Definition \ref{Def_deviation}, the function $I\ni t \mapsto  \sum_{i=1}^n D(x_i,t)$ is strictly decreasing and continuous.
Further, by (i) and (ii) of Definition \ref{Def_deviation}, if $\min(x_1,\ldots,x_n) = \max(x_1,\ldots,x_n)$, then
the unique solution $t\in I$ of the equation $\sum_{i=1}^n D(x_i,t) = n D(x_1,t) = 0$, is $t_0=x_1$.
On the other hand, if $\min(x_1,\ldots,x_n) < \max(x_1,\ldots,x_n)$, then
 \begin{align*}
    \sum_{i=1}^n D(x_i,t) \begin{cases}
                              > 0 \quad \text{if $t \leq  \min(x_1,\ldots,x_n)$,}\\
                              < 0 \quad \text{if $t \geq  \max(x_1,\ldots,x_n)$,}
                          \end{cases}
 \end{align*}
and hence, by the intermediate value theorem, there exists a unique $t_0$ in the open interval $(\min(x_1,\ldots,x_n),\max(x_1,\ldots,x_n))$
such that $\sum_{i=1}^n D(x_i,t_0)=0$.
So the notion of a deviation mean according to Definition \ref{Def_deviation_mean} is well-defined.
Further, for each $n\in\NN$, $\mathscr{M}_n^D$ is a strict, symmetric $n$-variable mean on $I$ in the sense of Definition \ref{Def_mean}, see Dar\'oczy \cite{Dar71b,Dar72b}.
In view of the results of Dar\'oczy and P\'ales \cite{DarPal82}, two deviation means on an open (nonempty) interval $I$ of $\RR$, generated by $D_1$ and $D_2$ are identical
 (i.e., for all $n\in\NN$, the equality $\mathscr{M}_n^{D_1}=\mathscr{M}_n^{D_2}$ holds on $I^n$) if and only if there exists a positive
 function $d:I\to\RR$ such that
\[
  D_1(x,t)=d(t)D_2(x,t), \qquad x,t\in I
\]
is valid. We note that this equality implies the continuity of $d$.
\proofend
\end{Rem}

In Remark \ref{Rem_M_estimate} we recall the notion of M-estimators and the more general concept of extremum estimators.
Then, in Remark \ref{Rem_M_estimate_deviation_mean}, we point out the fact that a deviation mean $\mathscr{M}_n^{D}(x_1,\ldots,x_n)$ of $x_1,\ldots,x_n$,
 where $x_1,\ldots,x_n$ are elements of an open (nonempty) interval can be considered as a unique solution to a convex optimization problem,
 yielding that deviation means are special M-estimators and $\psi$-estimators as well.
The optimization problem in question is somewhat similar to the one appearing in the definition of Fr\'echet means as well.

\begin{Rem}\label{Rem_M_estimate}
Let $(X,\cX)$ be a measurable space, $\Theta$ be a Borel subset of $\RR$, and $\varrho:X\times\Theta\to\RR$
 be a function such that for each $\theta\in\Theta$, the function $X\ni x\mapsto \varrho(x,\theta)$ is measurable
 with respect to the sigma-algebras $\cX$ and $\cB(\RR)$, where $\cB(\RR)$ denotes the Borel sigma-algebra on $\RR$.
Let $(\xi_n)_{n\in\NN}$ be a sequence of i.i.d.\ random variables with values in $X$ such that the distribution
 of $\xi_1$ depends on an unknown parameter $\theta\in\Theta$.
For each $n\in\NN$, Huber \cite{Hub64, Hub67} introduced a natural estimator of $\vartheta$ based on the observations
 $\xi_1,\ldots,\xi_n$ as a solution $\widehat\vartheta_n:=\widehat\vartheta_n(\xi_1,\ldots,\xi_n)$ of the following minimization problem:
 \begin{align}\label{help_M_est_min_problem}
   \argmin_{\vartheta\in\Theta}\frac{1}{n} \sum_{i=1}^n \varrho(\xi_i,\vartheta),
   \qquad \text{i.e.,}\qquad
   \inf_{\vartheta\in\Theta} \frac{1}{n} \sum_{i=1}^n \varrho(\xi_i,\vartheta)
     = \sum_{i=1}^n \varrho(\xi_i,\widehat\vartheta_n),
 \end{align}
 provided that such a solution exits.
One calls $\widehat\vartheta_n$ an M-estimator of the unknown parameter $\vartheta\in\Theta$ based on
 the i.i.d.\ observations $\xi_1,\ldots,\xi_n$, where M refers to ''maximum likelihood type''.
Under suitable regularity assumptions, the minimization problem in question can be solved by setting the derivative
 of the objective function (with respect to the unknown parameter) equal to zero:
 \[
 \frac{1}{n}\sum_{i=1}^n \partial_2\varrho(\xi_i,\vartheta)=0, \qquad \vartheta\in\Theta.
 \]
In the statistical literature, $\partial_2\varrho$ is often denoted by $\psi$, and hence in this case the M-estimator
 is often called $\psi$-estimator, while other authors call it a Z-estimator (the letter Z refers to ''zero'').

The class of M-estimators is in fact a special class of extremum estimators that we recall now.
For each $n\in\NN$, let $Q_n:X^n\times\Theta\to \RR$ be given functions, and
 an extremum estimator of $\theta\in\Theta$ based on the observations $\xi_1,\ldots,\xi_n$ is defined as
 a solution $\widehat\vartheta_n:=\widehat\vartheta_n(\xi_1,\ldots,\xi_n)$ of the following minimization problem:
 \[
   \argmin_{\theta\in\Theta} Q_n(\xi_1,\dots,\xi_n,\theta),
      \qquad \text{i.e.,}\qquad
    \inf_{\vartheta\in\Theta} Q_n(\xi_1,\dots,\xi_n,\theta)
      = Q_n(\xi_1,\dots,\xi_n,\widehat\vartheta_n),
 \]
 provided that such a solution exists.
\proofend
\end{Rem}

\begin{Rem}\label{Rem_M_estimate_deviation_mean}
We point out that a deviation mean on an open (nonempty interval) is a special M-estimator and a $\psi$-estimator as well in the language of statistics.
Let $I$ be an open (nonempty) interval of $\RR$, and let $D:I^2\to\RR$ be a deviation.
Using the notations of Remark \ref{Rem_M_estimate}, let $X:=I$, $\Theta:=I$, and define $\varrho:I^2\to\RR$,
 $$
   \varrho(x,t):=\int^t_x -D(x,s)\,\dd s, \qquad t,x\in I.
 $$
Then $\varrho(x,x)=0$, $x\in I$, and, by part (i) of Definition \ref{Def_deviation},
 for each $x\in I$, the function $t\ni I\mapsto \varrho(x,t)$ is nonnegative and differentiable and
 $$
   \partial_2 \varrho(x,t)=-D(x,t), \qquad t\in I.
 $$
Hence for each $x\in I$, the function $t\ni I\mapsto \varrho(x,t)$ is strictly convex,
 and therefore, since $\varrho(x,x)=0$, it is strictly decreasing on $I\cap(-\infty,x)$ and is strictly increasing on $I\cap(x,\infty)$.
Consequently, for each $n\in\NN$ and $x_1,\dots,x_n\in I$, there exists a unique minimum point $t_0$ of the strictly convex, differentiable function
 $$
  I \ni t\mapsto\varrho(x_1,t)+\dots+\varrho(x_n,t),
 $$
 and $t_0\in[\min(x_1,\dots,x_n),\max(x_1,\dots,x_n)]$.
Therefore, for each $n\in\NN$ and $x_1,\dots,x_n\in I$, we have
 $$
  \partial_2(\varrho(x_1,t)+\dots+\varrho(x_n,t))|_{t=t_0}=0,
 $$
 that is,
 $$
  D(x_1,t_0)+\dots+D(x_n,t_0)=0,
 $$
 and thus $t_0=\mathscr{M}_n^{D}(x_1,\ldots,x_n)$.
All in all, for each $n\in\NN$ and $x_1,\dots,x_n\in I$, the deviation mean $\mathscr{M}_n^{D}(x_1,\ldots,x_n)$ of $x_1,\ldots,x_n$ is nothing else
 but the unique solution of the minimization problem
 \[
    \min_{t\in I} \big(\varrho(x_1,t)+\dots+\varrho(x_n,t)\big).
 \]
In the particular case when $D(x,t)=x-t$, $x,t\in\RR$, the corresponding function $\varrho$ takes
 the form $\varrho(x,t)= \int_x^t (s-x)\,\dd s = \frac{1}{2}(x-t)^2$, $x,t\in\RR$,
 and hence, the above minimization problem leads to the least squares method.
\proofend
\end{Rem}

The class of deviation means is quite a broad class of means.
One can easily check that an $n$-variable Bajraktarevi\'c mean generated by the functions $f$ and $p$ (given in Definition \ref{Def_Bajraktarevic}) with $f$ being strictly increasing is an $n$-variable deviation mean corresponding to the deviation $D:I^2\to\RR$, $D(x,t):=p(x)(f(x)-f(t))$, $x,t\in I$.
P\'ales \cite{Pal82a,Pal85a,Pal88a,Pal91} studied several inequalities, algebraic and analytical properties, and obtained a characterization of (quasi)deviation means.
For a recent generalization and a survey on already existing generalizations of deviation means, see Stup\v{n}anov\'a and Smrek \cite{StuSmr}.
For an application of quasi-arithmetic and Bajraktarevi\'c means to congressional apportionment in the USA's election, see Appendix D in Barczy and Burai \cite{BarBur}.
Quasideviation means, which generalize the notion of deviation means, also have some applications in insurance mathematics and in the theory of lotteries.
Chudziak \cite{Chu} proved that the zero utility principle in insurance mathematics is a special case of quasideviation means,
 and, as a consequence, a new tool was demonstrated for dealing with the properties of this premium principle.
Chudziak and Chudziak \cite{ChuChu} showed that willingness to accept and willingness to pay in the theory of lotteries
 are particular cases of quasideviation means.

There is huge literature on limit theorems for M-estimators and extremum estimators based on i.i.d.\ observations
 and on more general observations as well, see, e.g., Huber \cite{Hub64, Hub67}, Huber and Ronchetti \cite{HubRon},
 He and Wang \cite{HeWan}, Rubin and Rukhin \cite{RubRuk}, Arcones \cite{Arc} and Hayashi \cite[Chapter 7]{Hay}.
In Remark \ref{Rem_feltet_osszehas}, we compare the conditions of our limit theorems with \emph{some} of the existing conditions in the literature.

Now we turn to recall some existing results on strong laws of large numbers and central limit theorems
 formulated for means of i.i.d.\ random variables, namely, in the setting of quasi-arithmetic and Bajraktarevi\'c means of i.i.d.\ random variables.
For quasi-arithmetic means of i.i.d.\ random variables, a strong law of large numbers readily follows by Kolmogorov's strong law of
 large numbers, and de Carvalho \cite[Theorem 1]{Car} derived a central limit theorem as well.
Before presenting de Carvalho's results, let us recall that if $f:I\to\RR$ is a continuous and strictly increasing function, where $I$ is a nondegenerate interval of $\RR$, and $\xi$ is a random variable such that $\PP(\xi\in I)=1$ and
 $\EE(\vert f(\xi)\vert)<\infty$, then  Kolmogorov's \emph{quasi-arithmetic expected value of $\xi$ corresponding to $f$} is defined by
 \[
   \mathbb{A}^f(\xi):=f^{-1} (\EE (f(\xi)) ).
 \]
Here $\EE (f(\xi))\in f(I)$ (for a proof, see Lemma \ref{Lem_Exp_conv_hull}),
 so the notion $\mathbb{A}^f(\xi)$ is well-defined.
Observe that if, for some $n\in\NN$, the random variable $\xi$ takes the values $x_1,\dots,x_n\in I$ with equal $1/n$ probabilities, then $\mathbb{A}^f(\xi)=\mathscr{A}_n^f(x_1,\dots,x_n)$, which motivates the terminology here.
If $I=(0,\infty)$ and $f(x) = x^p$, $x>0$, where $p\geq 1$, and $\EE(\xi^p)<\infty$, then $\mathbb{A}^f(\xi) = (\EE(\xi^p))^{\frac{1}{p}}$, which is nothing else, but the $L_p$-norm of $\xi$. The standard expected value of $\xi$ corresponds to $f:\RR\to\RR$ given by $f(x):=ax+b$, $x\in\RR$, where $a,b\in\RR$, $a\ne 0$.
Recall also that $\var(\xi):= \EE((\xi - \EE(\xi))^2)$ whenever $\EE(\vert\xi\vert)<\infty$.

\begin{Thm}[de Carvalho (2016)]\label{Thm_de_Carvalho}
Let $I$ be a nondegenerate  interval of $\RR$, and $f:I\to\RR$ be a continuous and strictly increasing function.
Let $(\xi_n)_{n\in\NN}$ be a sequence of i.i.d.\ random variables such that
 $\PP(\xi_1\in I)=1$ and $\EE(|f(\xi_1)|)<\infty$.
Then
 \[
  \mathscr{A}_n^f(\xi_1,\ldots,\xi_n) \as \mathbb{A}^f(\xi_1) \qquad \text{as $n\to\infty$,}
 \]
If, in addition, $ \var( f(\xi_1) )\in(0,\infty)$ and $f'(\mathbb{A}^f(\xi_1))$  exists and is non-zero, then
 \[
   \sqrt{n}\big(\mathscr{A}_n^f(\xi_1,\ldots,\xi_n) - \mathbb{A}^f(\xi_1)\big) \distr \cN\left(0, \frac{\var(f(\xi_1))}{( f'(\mathbb{A}^f(\xi_1)))^2} \right)
   \qquad \text{as $n\to\infty$.}
 \]
\end{Thm}

Very recently, Barczy and Burai \cite[Theorem 2.1]{BarBur} have proved a strong law of large numbers and a central limit theorem for the Bajraktarevi\'c means of i.i.d.\ random variables.

\begin{Thm}[Barczy and Burai (2022)]\label{Thm_CLT_Baj_mean}
Let $I$ be a nondegenerate interval of $\RR$, let $f:I\to\RR$ be a continuous and strictly monotone
 function such that  the interval $f(I)$ is closed, and let $p:I\to (0,\infty)$ be a Borel measurable (weight) function.
Let $(\xi_n)_{n\in\NN}$ be a sequence of i.i.d.\ random variables such that $\PP(\xi_1\in I)=1$, $\EE(p(\xi_1))<\infty$ and $\EE(p(\xi_1)|f(\xi_1)|)<\infty$. Then
 \[
    \mathscr{B}^{f,p}_n(\xi_1,\ldots,\xi_n) \as f^{-1}\left( \frac{\EE(p(\xi_1)f(\xi_1))}{\EE(p(\xi_1))} \right)
    =:\mathbb{B}^{f,p}(\xi_1)\qquad \text{as $n\to\infty$.}
 \]
If, in addition, $\EE((p(\xi_1))^2)<\infty$, $\EE((p(\xi_1) f(\xi_1))^2)<\infty$, and $f$ is differentiable at $\mathbb{B}^{f,p}(\xi_1)$ with a non-zero derivative, then
 \begin{align}\label{help7_Baj_mean}
  \sqrt{n} \left(\mathscr{B}^{f,p}_n(\xi_1,\ldots,\xi_n) - \mathbb{B}^{f,p}(\xi_1) \right) \distr \cN\big(0, \sigma^2_{f,p}\big)
    \qquad \text{as $n\to\infty$,}
 \end{align}
where
 \begin{align*}
  \sigma_{f,p}^2
    := \frac{(\EE(p(\xi_1)))^{-4}}
            {\left(f'(\mathbb{B}^{f,p}(\xi_1))\right)^2}
         &\Big((\EE(p(\xi_1)))^2 \var(p(\xi_1) f(\xi_1)) + (\EE(p(\xi_1) f(\xi_1) ))^2 \var(p(\xi_1)) \\
         &\phantom{\Big(\;} - 2 \EE(p(\xi_1)) \EE(p(\xi_1)f(\xi_1)) \cov(p(\xi_1), p(\xi_1)f(\xi_1)) \Big).
 \end{align*}
\end{Thm}

Here we call the attention to the fact that in the above recalled Theorem \ref{Thm_CLT_Baj_mean} of Barczy and Burai \cite[Theorem 2.1]{BarBur}
 it is supposed that $f(I)$ is closed, however, as a consequence of our forthcoming Theorems \ref{Thm_SLLN_deviation} and \ref{Thm_CLT_deviation},
 it turns out that is a superfluous assumption.

Recently, there has also been a renewed interest in operator means in functional analysis, here we only note that Lim and P\'alfia \cite{LimPal} proved a strong law of large numbers for the $L^1$-Karcher mean.

In this paper we establish a strong law of large numbers, a central limit theorem, a law of the iterated logarithm
and large deviation theorem for deviation means of i.i.d.\ random variables without referring to general limit theorems for M-estimators and extremum estimators, see Theorems \ref{Thm_SLLN_deviation}, \ref{Thm_CLT_deviation}, \ref{Thm_LIL_deviation} and \ref{Thm_Dev_mean_large_dev}, respectively
 (for the strong law of large numbers, we suppose only pairwise independence instead of (total) independence).
Our results on the strong law of large numbers and on the central limit theorem generalize the corresponding ones for quasi-arithmetic means
 due to de Carvalho \cite[Theorem 1]{Car} and the ones for Bajraktarevi\'c means due to Barczy and Burai \cite[Theorem 2.1]{BarBur}.
The paper is organized as follows.
In Section \ref{section_dev_mean_rv}, we introduce the notion of deviation mean of a random variable and we give several examples as well.
Sections \ref{section_SLLN}, \ref{section_CLT}, \ref{section_LIL} and \ref{section_LD} contain our limit theorems and, in Section \ref{section_ELT}, we demonstrate their specializations to Bajraktarevi\'c means and to a special deviation mean which is not a Bajraktarevi\'c mean (see Examples \ref{Ex_3_Bajraktarevic} and \ref{Ex_4_nonBajraktarevic}).
The proof of the strong law of large numbers for deviation means  can be traced back to Etemadi's strong law of large numbers for pairwise independent and identically distributed random variables. For proving a central limit theorem and a law of the iterated logarithm for deviation means, as an auxiliary result, we derive some sufficient conditions on a deviation $D$ and on a sequence of i.i.d.\ random variables $(\xi_n)_{n\in\NN}$ under which
 \[
  \tau_n \sum_{i=1}^n \Big(D\big(\xi_i, t_0 + \kappa_n \big)- D(\xi_i, t_0)
                                        - \kappa_n \partial_2D(\xi_i, t_0)\Big)
       \as 0 \qquad \text{as \ $n\to\infty$}
 \]
holds with some appropriate sequences $(\kappa_n)_{n\in\NN}$ and $(\tau_n)_{n\in\NN}$ of real numbers, where $t_0$ is an interior point of the nondegenerate interval $I$, see Lemma \ref{Lem_aux_dev_as}.
The weak convergence in the central limit theorem for deviation means is proved by using the definition of weak convergence together with classical central limit theorem for i.i.d.\ random variables, the strong law of large numbers due to Marcinkiewicz and Zygmund, Slutsky's lemma and the above described Lemma \ref{Lem_aux_dev_as}.
The proof of the law of the iterated logarithm in terms of deviation means is, in part, based on the classical law of the iterated logarithm for i.i.d.\ random variables and on Lemma~\ref{Lem_aux_dev_as}.
The proof of the large deviation theorem for deviation means can be easily traced back to Cram\'er's theorem on large deviations.
We close the paper with an appendix, where, for the convenience of the reader, we recall or prove some auxiliary results that are used in the proofs such as the strong law of large numbers due to Marcinkiewicz and Zygmund, Lemma \ref{Lem_aux_dev_as} and Cram\'er's theorem on large deviations.

\section{Deviation mean of a random variable}
\label{section_dev_mean_rv}

First we introduce the deviation mean of a random variable.

\begin{Def}
Let $I\subset\RR$ be a nondegenerate interval and $D:I^2\to\RR$ be a deviation. We say that $D$ is Borel measurable in the first variable if, for all $t\in I$, the function $I\ni x\mapsto D(x,t)$ is Borel measurable.
\end{Def}

\begin{Lem}\label{Lem_Dev_mean_Borel_measurable}
Let $I\subset\RR$ be a nondegenerate interval and $D:I^2\to\RR$ be a deviation which is Borel measurable in the first variable.
Then, for all $n\in\NN$, the $n$-variable $D$-deviation mean $\mathscr{M}_n^D: I^n\to I$ is Borel measurable. Furthermore, if $n\in\NN$ and $\xi_1,\dots,\xi_n$ are random variables on a probability space $(\Omega,\cF,\PP)$ with values in $I$, then $\mathscr{M}_n^D(\xi_1,\dots,\xi_n)$ is a random variable, that is, the map
\begin{equation}\label{eq:Mxi}
  \Omega\ni\omega\mapsto\mathscr{M}_n^D(\xi_1(\omega),\dots,\xi_n(\omega))
\end{equation}
is $\cF$-measurable.
\end{Lem}

\noindent{\bf Proof.}
Let $t\in I$ be arbitrary. Then the function
\[
  I^n\ni(x_1,\dots,x_n)\mapsto \sum_{i=1}^n D(x_i,t)
\]
is Borel measurable, and therefore the level set
\[
 \big\{(x_1,\dots,x_n)\in I^n\colon\mathscr{M}_n^D(x_1,\dots ,x_n)<t\big\}
 =\bigg\{(x_1,\dots, x_n)\in I^n\colon\sum_{i=1}^n D(x_i,t)<0\bigg\}
\]
is also Borel measurable. In order to establish this equality, for each $x_1,\ldots,x_n\in I$, observe that the function $I\ni t \mapsto \sum_{i=1}^n D(x_i,t)$ is strictly decreasing and vanishes at the point $t=\mathscr{M}_n^D(x_1,\dots ,x_n)$.

To see that the second assertion is also valid, observe that the map given by \eqref{eq:Mxi} is the composition of a Borel measurable and of an $\cF$-measuarable function, and hence, indeed, it is a random variable.
\proofend

\begin{Lem}\label{Lem_Dev_mean_random_1}
Let $I\subset\RR$ be a nondegenerate interval and $D:I^2\to\RR$ be a deviation which is Borel measurable in the first variable.
 Let $\xi$ be a random variable with $\PP(\xi\in I)=1$.
Let
 \[
 I_{\xi}^D:=\{t\in I\mid \EE(|D(\xi,t)|)<\infty\}.
 \]
Then $I_{\xi}^D$ is a (possibly degenerate) subinterval of $I$ and the function $g:I_{\xi}^D\to\RR$ defined by $g(t):=\EE(D(\xi,t))$, $t\in I_{\xi}^D$, is continuous and strictly decreasing.
Consequently, if there exists $t_0\in I_{\xi}^D$ with $g(t_0)=0$, then $t_0$ is unique.
\end{Lem}

\noindent{\bf Proof.} Let $a,b\in I_{\xi}^D$ with $a<b$ and let $t\in[a,b]$.
Then, by the monotonicity property of deviations (see part (i) of Definition \ref{Def_deviation}), we have
\[
  D(\xi,a)\geq D(\xi,t)\geq D(\xi,b).
\]
Therefore,
\[
  |D(\xi,t)|\leq\max(|D(\xi,a)|,|D(\xi,b)|)
  \leq |D(\xi,a)|+|D(\xi,b)|,
\]
which shows that $\EE(|D(\xi,t)|)<\infty$. Hence $t\in I_{\xi}^D$, proving that $I_{\xi}^D$ is an interval.

By Lebesgue's Dominated Convergence Theorem, it immediately follows that the function $g$ is continuous.
If $s,t\in I_{\xi}^D$ with $s < t$, then, by the monotonicity property of deviations, we have $D(x,s)> D(x,t)$, $x\in I$, yielding that $D(\xi(\omega),s)> D(\xi(\omega),t)$ for almost every $\omega\in\Omega$.
Hence $\EE(D(\xi,s))\geq \EE(D(\xi,t))$, and $\EE(D(\xi,s)) = \EE(D(\xi,t))$ cannot hold, otherwise the expectation of the almost surely positive random variable $D(\xi,s)-D(\xi,t)$ would be zero. Thus $g$ is strictly decreasing.
\proofend

\begin{Ex}\label{Ex_1}
(i). If $\xi$ takes only finitely many values of $I$, then $\EE(\vert D(\xi,t)\vert)<\infty$ holds for $t\in I$, and hence $I^D_\xi = I$.
Note also that in this special case we do not need to assume that $D$ is Borel measurable in the first variable in order that
 $D(\xi,t)$ be a random variable for all $t\in I$.

\noindent (ii). Let $I\subset\RR$ be a nondegenerate interval and let $D:I^2\to\RR$, $D(x,t):=p(x)(f(x)-f(t))$, $x,t\in I$,
 where $f:I\to\RR$ is a continuous and strictly increasing function and $p:I\to (0,\infty)$ is a Borel measurable function.
Then $D$ is deviation which is Borel measurable in the first variable.
Further, let $\xi$ be a random variable such that $\PP(\xi\in I)=1$, $\EE(p(\xi))<\infty$ and $\EE(p(\xi)\vert f(\xi)\vert )<\infty$.
Then $\EE(\vert D(\xi,t)\vert)<\infty$ for each $t\in I$, so $I^D_\xi = I$. In fact, it is not difficult to see that the following reversed implication is also valid: if $\EE(\vert D(\xi,t)\vert)<\infty$ holds for at least two distinct values of $t\in I$, then $\EE(p(\xi))<\infty$ and $\EE(p(\xi)\vert f(\xi)\vert )<\infty$.

\noindent (iii). Let $I\subset\RR$ be a compact interval, $D:I^2\to\RR$ be a continuous deviation and $\xi$ be a random variable such that $\PP(\xi\in I)=1$. Then $\EE(\vert D(\xi,t)\vert)<\infty$, $t\in I$, since $D$ is bounded on the compact set $I^2$, and hence $I^D_\xi = I$.

\noindent (iv). Let $I:=\RR$ and $D:I^2\to\RR$, $D(x,t):=x-t$, $x,t\in I$.
Further, let $\xi$ be a random variable such that $\EE(\vert \xi\vert) = \infty$.
Then $\EE(\vert D(\xi,t)\vert) = \EE(\vert \xi-t\vert) = \infty$, $t\in I$, and thus $I^D_\xi=\emptyset$.

\noindent (v). It is not too difficult to give an example for a nondegenerate interval $I$,
 for a deviation $D$ which is Borel measurable in the first variable
 and for a random variable $\xi$ such that $I^D_\xi$ is nondegenerate and does not coincide with $I$. Let $I:=[-2,\infty)$, and $D:I^2\to \RR$,
 \[
   D(x,t):=\begin{cases}
              \ee^{xt} - \ee^{x^2} & \text{if \ $x\in[{-2},0)$ \ and \ $t\geq {-2}$,}\\
              x-t=-t & \text{if \ $x=0$ \ and \ $t\geq {-2}$,}\\
              \ee^{-xt} - \ee^{-x^2} & \text{if \ $x>0$ \ and \ $t\geq {-2}$.}
           \end{cases}
 \]
Then $D$ is a deviation which is Borel measurable in the first variable, since
 \[
  \partial_2 D(x,t) = \begin{cases}
                        x\ee^{tx} <0& \text{if \ $x\in[{-2},0)$ \ and \ $t\geq {-2}$,}\\
                        -1 <0& \text{if \ $x=0$ \ and \ $t\geq {-2}$,}\\
                        -x\ee^{-xt} < 0& \text{if \ $x>0$ \ and \ $t\geq {-2}$.}
                     \end{cases}
 \]
Let $\xi:=\ee^{\eta}$, where $\eta$ is a standard normally distributed random variable, i.e., $\xi$ has a lognormal distribution with parameters $0$ and $1$ and it has a density function
 \[
   f_\xi(x)
   :=\begin{cases}
       \frac{1}{\sqrt{2\pi}x} \ee^{-\frac{1}{2}(\ln(x))^2} & \text{if \ $x>0$,}\\
       0 & \text{if \ $x\leq 0$.}
     \end{cases}
 \]
Then $\PP(\xi\in (0,\infty))=1$, and hence \ $\PP(\xi\in I)=1$ holds as well.
Further, using that $\PP(\xi\in (0,\infty))=1$, we have
 \begin{align*}
   \EE(\vert D(\xi,t)\vert)
     = \EE( \vert \ee^{-\xi t} - \ee^{-\xi^2}\vert), \qquad t\in I.
 \end{align*}
If $t\geq 0$, then $\EE( \vert \ee^{-\xi t} - \ee^{-\xi^2}\vert) \leq 1+1=2<\infty$.
If $t\in[-2,0)$, then, using that $\vert u-v\vert\geq \vert \vert u\vert - \vert v\vert \vert$, $u,v\in\RR$, we have
 \begin{align*}
   \EE( \vert \ee^{-\xi t} - \ee^{-\xi^2}\vert)
     \geq \EE\big(\ee^{-\xi t} - \ee^{-\xi^2} \big)
     = \EE(\ee^{-\xi t}) - \EE(\ee^{-\xi^2}),
 \end{align*}
 where $\EE(\ee^{-\xi^2})\in(0,1)$, and
 \[
   \EE(\ee^{-\xi t}) = \int_0^\infty \ee^{-xt} \frac{1}{\sqrt{2\pi}x} \ee^{-\frac{1}{2}(\ln(x))^2}\,\dd x
                     =  \frac{1}{\sqrt{2\pi}} \int_0^\infty  \frac{1}{x} \ee^{-xt - \frac{1}{2}(\ln(x))^2}\,\dd x.
 \]
Using L'Hospital's rule one can check that if $t\in[-2,0)$, then
 \[
   \lim_{x\to\infty} \left( -xt - \frac{1}{2}(\ln(x))^2 \right)   = \infty,
 \]
 and consequently, $\lim_{x\to\infty} \ee^{-xt - \frac{1}{2}(\ln(x))^2} = \infty$.
Hence for each $L>0$ there exists $x_L>0$ (may depend on $t$) such that $\ee^{-xt - \frac{1}{2}(\ln(x))^2} > L$ if $x\geq x_L$.
Thus for each $L>0$, we have
 \[
    \int_0^\infty  \frac{1}{x} \ee^{-xt - \frac{1}{2}(\ln(x))^2}\,\dd x
      \geq L \int_{x_L}^\infty \frac{1}{x}\,\dd x = L\Big(\lim_{x\to\infty} \ln(x) - \ln(x_L)\Big) =\infty,
 \]
 yielding $\EE(\ee^{-\xi t})=\infty$ for $t\in[-2,0)$ and then $\EE( \vert \ee^{-\xi t} - \ee^{-\xi^2}\vert ) =\infty$ for $t\in[-2,0)$.
Therefore $I^D_\xi = [0,\infty)$.
Note that in this example $I^D_\xi$ is nondegenerate and does not coincide with $I=[-2,\infty)$.
\proofend
\end{Ex}

\begin{Def}\label{Def_deviation_mean_random}
Let $I\subset\RR$ be a nondegenerate interval and $D:I^2\to\RR$ be a deviation which is Borel measurable in the first variable.
Let $\xi$ be a random variable with $\PP(\xi\in I)=1$ and assume that there exists a point $t_0\in I_{\xi}^D$ such that $\EE(D(\xi,t_0))=0$.
This value $t_0$ is called the \emph{$D$-deviation mean of $\xi$} and is denoted by $\mathbb{M}^D(\xi)$. Furthermore, if $\mathbb{M}^D(\xi)$ belongs to the interior of the interval $I_{\xi}^D$, then we say that \emph{$\xi$ satisfies the interior point condition with respect to $D$}.
\end{Def}

\begin{Rem}
Under the conditions of Definition \ref{Def_deviation_mean_random}, by Lemma \ref{Lem_Dev_mean_random_1}, we have that if $\mathbb{M}^D(\xi)$ exists,
 then it is uniquely determined. Clearly, $\mathbb{M}^D(\xi)$ is determined by $D$ and the distribution of $\xi$.
Further, in case of an open (nonempty) interval $I$, the $D$-deviation mean $t_0=\mathbb{M}^D(\xi)$ of $\xi$ can be considered
 as the unique minimizer of the strictly convex and differentiable function
 $I\ni t \mapsto \int_{t_0}^t \EE(-D(\xi,s))\,\dd s$.
\proofend
\end{Rem}

\begin{Ex}\label{Ex_2}
\noindent (i). In the particular case when $\xi$ takes the values $x_1,\dots,x_n\in I$ with equal $\frac1n$ probabilities, where $n\in\NN$, one can check that
 \[
  \mathbb{M}^D(\xi)=\mathscr{M}_n^{D}(x_1,\ldots,x_n).
 \]

\noindent (ii).
Let $p>0$, $I:=\RR$, and $D:I^2\to\RR$ be given by $D(x,t):=\sign(x-t)\cdot |x-t|^p$, $x,t\in I$.
Then $D$ is a deviation which is Borel measurable in the first variable.
Further, for any fixed $x\in I$, if $p=1$, then $D(x,t)=x-t$ and hence $\partial_2 D(x,x)=-1$; if $p>1$, then $\partial_2 D(x,x)=0$;
 and if $p<1$, then $\partial_2 D(x,x)$ does not exist.
If $\xi$ is a random variable and $t\in I$, then $\EE(\vert D(\xi,t)\vert)<\infty$ holds if and only if
 $\EE(\vert \xi-t\vert^p)<\infty$, which is equivalent to $\EE(\vert \xi\vert^p)<\infty$.
Consequently,
 \begin{align}\label{help_example1}
     I_\xi^D = \begin{cases}
                  I & \text{if \ $\EE(\vert \xi\vert^p)<\infty$},\\
                   \emptyset & \text{if \ $\EE(\vert \xi\vert^p)=\infty$.}
                \end{cases}
 \end{align}

Now, consider the case when $p:=2$, and $x_1,\dots,x_n\in I$ with $x_1\leq \dots\leq x_n$ and $x_1<x_n$, where $n\in\NN\setminus\{1\}$.
Let $s:=\mathscr{M}_n^D(x_1,\dots,x_n)$. Then, for some $j\in\{1,\dots,n-1\}$, we have that $x_j\leq s<x_{j+1}$, therefore
\[
   -\sum_{i=1}^j (x_i-s)^2+\sum_{i=j+1}^n (x_i-s)^2=0,
\]
equivalently,
\[
   (n-2j)s^2+2\bigg(\sum_{i=1}^j x_i-\sum_{i=j+1}^n x_i\bigg)s
   +\bigg(\sum_{i=j+1}^n x_i^2-\sum_{i=1}^j x_i^2\bigg)=0.
\]
If $n=2$, then $x_1<x_2$, $j=1$, and $s$ equals the arithmetic mean of $x_1$ and $x_2$.
For $n\geq3$, $n\in\NN$,  the mean $\mathscr{M}_n^D(x_1,\dots,x_n)$ is the unique solution $s$ of the above equation.

Now let $p:=2$ and $\xi$ be an exponentially distributed random variable with parameter $1$.
Then $\EE(\xi^k)=k!$, $k\in\ZZ_+$, and hence, by \eqref{help_example1}, $I_\xi^D=I$.
Moreover, if $t\leq 0$, then
  \begin{align*}
     \EE(D(\xi,t))
       & = \EE(\sign(\xi-t)\cdot(\xi-t)^2)
         = \int_0^\infty (x-t)^2 \ee^{-x}\,\dd x
         = \int_0^\infty (x^2-2tx+t^2) \ee^{-x}\,\dd x \\
       & = \EE(\xi^2) - 2t\EE(\xi) + t^2
         = 2 -2t+t^2
         = (t-1)^2 +1,
  \end{align*}
 and if $t>0$, then
 \begin{align*}
     \EE(D(\xi,t))
       & = \EE(\sign(\xi-t) \cdot(\xi-t)^2)
         = - \int_0^t (x-t)^2 \ee^{-x}\,\dd x  + \int_t^\infty (x-t)^2 \ee^{-x}\,\dd x\\
       & = - \int_0^\infty (x-t)^2 \ee^{-x}\,\dd x + 2\int_t^\infty (x-t)^2 \ee^{-x}\,\dd x\\
       & = -(t-1)^2 - 1 + 2\ee^{-t}\int_t^\infty (x-t)^2 \ee^{-(x-t)}\,\dd x \\
       & = -(t-1)^2 - 1 + 2\ee^{-t}\int_0^\infty y^2 \ee^{-y}\,\dd y \\
       & = -(t-1)^2 - 1 + 2\ee^{-t}\EE(\xi^2)
        = 4\ee^{-t} - (1-t)^2 -1.
  \end{align*}
Consequently, the equation $\EE(D(\xi,t))=0$, $t\in I$, takes the form
 \[
   4\ee^{-t} = (1-t)^2 +1, \qquad t>0.
 \]
This equation has a unique (real) solution (denoted by $\mathbb{M}^D(\xi)$),
 and it is approximately $1.300075$ (using software R).

Now let $p:=2$ and $\xi$ be an absolutely continuous random variable having a density function
 \[
   f_\xi(x):=\begin{cases}
               3x^{-4} & \text{if \ $x>1$,}\\
               0 & \text{if \ $x\leq 1$.}
             \end{cases}
 \]
Then $\EE(\xi^2) = 3\int_1^\infty x^{-2}\,\dd x= 3<\infty$, and hence, by \eqref{help_example1}, $I_D^\xi = I$.
Further, if $t\leq 1$, then
  \begin{align*}
     \EE(D(\xi,t))
       & = \EE(\sign(\xi-t)\cdot(\xi-t)^2)
         = 3\int_1^\infty (x-t)^2 x^{-4}\,\dd x  \\
       & = 3\int_1^\infty x^{-2}\,\dd x - 6t\int_1^\infty x^{-3}\,\dd x
           + 3t^2 \int_1^\infty x^{-4}\,\dd x
         =  3\Big(1 - t + \frac{t^2}{3}\Big),
  \end{align*}
 and if $t>1$, then we have
 \begin{align*}
     \EE(D(\xi,t))
       & = \EE(\sign(\xi-t)\cdot(\xi-t)^2)
         =  \int_1^t (-1)(x-t)^2 3x^{-4}\,\dd x  + \int_t^\infty (x-t)^2 3x^{-4}\,\dd x\\
       & = -3 \int_1^\infty (x-t)^2 x^{-4}\,\dd x
           + 6\int_t^\infty (x-t)^2 x^{-4}\,\dd x \\
       &= -3\Big(1 - t + \frac{t^2}{3}\Big) + \frac{2}{t} .
  \end{align*}
Consequently, since $\EE(D(\xi,1))=1>0$, the equation $\EE(D(\xi,t))=0$, $t\in I$, takes the form
 \[
    -t^3 + 3t^2 - 3t +2 =0,\qquad t\geq 1,
 \]
 which is equivalent to
 \[
   (t-2)(t^2-t+1)=0,\qquad t\geq 1.
 \]
This equation has a unique (real) solution $2$, that is, $\mathbb{M}^D(\xi)=2$.

\noindent (iii).
We give an example for a nondegenerate interval $I$, for a deviation $D$ and for a random random variable $\xi$ such that
 $I^D_\xi$ is nondegenerate, $I^D_\xi$ does not coincide with $I$ and there does not exist a $t_0\in I^D_\xi$ such that
 $\EE(D(\xi,t_0))=0$.
Let $I:=[-2,\infty)$ and $D:I^2\to\RR$ be the deviation given in part (v) of Example \ref{Ex_1}.
As we checked, $D$ is a deviation which is Borel measurable in the first variable.
Let $\zeta:=\ee^{\eta}$, where $\eta$ is a standard normally distributed random variable,
 and let $\xi:=\zeta-2$.
Then $\PP(\xi\in(-2,\infty))=1$ and hence $\PP(\xi\in I)=1$ holds as well.
Further, for each $t\in I$, we have
 \begin{align*}
  \EE(\vert D(\xi,t)\vert)
   & = \EE(\vert D(\xi,t)\vert \bone_{\{\xi\leq 0\}}) + \EE(\vert D(\xi,t)\vert \bone_{\{\xi>0\}})\\
   & = \EE(\vert \ee^{\xi t} - \ee^{\xi^2} \vert \bone_{\{\xi\leq 0\}}) + \EE(\vert \ee^{-\xi t} - \ee^{-\xi^2} \vert \bone_{\{\xi>0\}}).
 \end{align*}
If $t\geq 0$, then
 \[
   \EE(\vert \ee^{\xi t} - \ee^{\xi^2} \vert \bone_{\{\xi\leq 0\}})
   \leq 1 + \ee^4<\infty,
 \]
 and $\EE(\vert \ee^{-\xi t} - \ee^{-\xi^2} \vert \bone_{\{\xi>0\}})\leq 1+1=2$,
 yielding that $\EE(\vert D(\xi,t)\vert)<\infty$, i.e., $t\in I^D_\xi$.
If $t\in[-2,0)$, then, using that $\vert u-v\vert\geq \vert \vert u\vert - \vert v\vert \vert$, $u,v\in\RR$, we have
 \begin{align*}
   \EE(\vert \ee^{-\xi t} - \ee^{-\xi^2} \vert \bone_{\{\xi>0\}})
     \geq \EE\big(\ee^{-\xi t}\bone_{\{\xi>0\}}\big)   - \EE\big(\ee^{-\xi^2} \bone_{\{\xi>0\}} \big),
 \end{align*}
 where, similarly as in part (v) of Example \ref{Ex_1} we have
 \[
   \EE\big(\ee^{-\xi t}\bone_{\{\xi>0\}}\big)
            = \ee^{2t}\EE(\ee^{-\zeta t} \bone_{\{\zeta>2\}})
            =\ee^{2t} \int_2^\infty \ee^{-xt} \frac{1}{\sqrt{2\pi}x} \ee^{-\frac{1}{2}(\ln(x))^2}\,\dd x
            = \infty,
 \]
 and $\EE\big(\ee^{-\xi^2} \bone_{\{\xi>0\}} \big)\in(0,1)$.
Consequently, if $t\in[-2,0)$, then $t\notin I^D_\xi$.
Hence $I^D_\xi=[0,\infty)$.

Moreover, by the dominated convergence theorem, we get
 \begin{align*}
  \lim_{t\to\infty} g(t)
   = \lim_{t\to\infty} \EE(D(\xi,t))
     & = \lim_{t\to\infty} \Big[\EE( (\ee^{\xi t} - \ee^{\xi^2} ) \bone_{\{\xi\leq 0\}}) + \EE( (\ee^{-\xi t} - \ee^{-\xi^2}) \bone_{\{\xi>0\}})\Big]\\
     & = - \EE( \ee^{\xi^2} \bone_{\{\xi\leq 0\}}) - \EE( \ee^{-\xi^2} \bone_{\{\xi>0\}})<0,
 \end{align*}
 and
 \begin{align*}
  g(0)
   &= \EE(D(\xi,0))
   = 1 - \EE( \ee^{\xi^2} \bone_{\{\xi\leq 0\}}) - \EE( \ee^{-\xi^2} \bone_{\{\xi>0\}})\\
   &= 1 - \int_0^2 \ee^{ (x-2)^2 } \frac{1}{\sqrt{2\pi} x}{\ee^{-\frac{1}{2}(\ln(x))^2}}\,\dd x
         - \int_2^\infty \ee^{ -(x-2)^2 } \frac{1}{\sqrt{2\pi} x}{\ee^{-\frac{1}{2}(\ln(x))^2}}\,\dd x
   <0.
 \end{align*}
Indeed, $\int_0^2 \ee^{ (x-2)^2 } \frac{1}{\sqrt{2\pi} x}{\ee^{-\frac{1}{2}(\ln(x))^2}}\,\dd x\approx 6,55323$
 and $\int_2^\infty \ee^{ -(x-2)^2 } \frac{1}{\sqrt{2\pi} x}{\ee^{-\frac{1}{2}(\ln(x))^2}}\,\dd x\approx 0,09372$.
By Lemma \ref{Lem_Dev_mean_random_1}, the function $g$ is continuous and strictly decreasing, and hence there does not exist $t_0\in I^D_\xi$ such that $g(t_0)=0$.

\noindent (iv).
Let $I:=[0,1]$ and $D:I^2\to\RR$, $D(x,t):=x-t$, $x,t\in I$.
Further, let $\xi$ be a random variable such that $\PP(\xi=0)=1$.
Then $\EE(\vert D(\xi,t)\vert) = \vert t\vert<\infty$, $t\in I$, and thus $I^D_\xi=I$.
Moreover, $\EE(D(\xi,t))=0$, $t\in I$, holds if and only if $t=0$.
Thus $\mathbb{M}^D(\xi)=0$ is a boundary point of $I_D^\xi$, and hence $\xi$ does not satisfy the interior point condition with respect to $D$.
\proofend
\end{Ex}

\section{Strong law of large numbers for deviation means}
\label{section_SLLN}

In the theorem below, we present a strong law of large numbers for deviation means of pairwise i.i.d.\ random variables.

\begin{Thm}\label{Thm_SLLN_deviation}
Let $I\subset\RR$ be a nondegenerate interval, let $D:I^2\to\RR$ be a deviation and $(\xi_n)_{n\in\NN}$ be a sequence of pairwise independent and identically distributed random variables.
Assume that
\vspace{-12pt}\begin{enumerate}[(i)]\itemsep=-4pt
 \item $D$ is Borel measurable in the first variable;
 \item $\PP(\xi_1\in I)=1$ and $\xi_1$ satisfies the interior point condition with respect to $D$.
\end{enumerate}\vspace{-12pt}
Then
 \begin{align}\label{help_dev_SLLN_allitas}
    \mathscr{M}_n^{D}(\xi_1,\ldots,\xi_n) \as \mathbb{M}^D(\xi_1) \qquad \text{as $n\to\infty$.}
 \end{align}
\end{Thm}

\noindent{\bf Proof.}
Since $\PP(\xi_1\in I)=1$, we have $\PP(\xi_i\in I, i\in\NN)=1$, and hence, without loss of generality, we may assume that
 $\xi_i(\omega)\in I$ for each $\omega\in \Omega$ and $i\in\NN$.
For each $n\in\NN$ and $\omega\in\Omega$, define $g_{n,\omega}: I \to\RR$ by
 \[
     g_{n,\omega}(t):=\frac{1}{n}\sum_{i=1}^n D(\xi_i(\omega),t),\qquad t\in I.
 \]
Further, let $g:I^D_{\xi_1}\to \RR$,
 \[
     g(t):=\EE(D(\xi_1,t)),\qquad t\in I^D_{\xi_1}.
 \]

By the properties of $D$ and Remark \ref{Rem1}, for each $n\in\NN$ and $\omega\in\Omega$, $g_{n,\omega}$ is continuous and strictly decreasing,
 and its unique zero equals $\mathscr{M}_n^{D}(\xi_1(\omega),\ldots,\xi_n(\omega))$.
Note that the random variables $D(\xi_n,t)$, $n\in\NN$, $t\in I$, are indeed well-defined (due to assumption (i) and $\PP(\xi_1\in I)=1$) and,
 for $t\in I^D_{\xi_1}$, they have finite first moments.
Further, due to Lemma \ref{Lem_Dev_mean_Borel_measurable}, for each $n\in\NN$, $\mathscr{M}_n^{D}(\xi_1,\ldots,\xi_n)$ is also a random variable.

In view of Lemma~\ref{Lem_Dev_mean_random_1}, $g$ is continuous and strictly decreasing on $I^D_{\xi_1}$. By the interior point condition, $t_0:=\mathbb{M}^D(\xi_1)$ exists and belongs to the interior of $I^D_{\xi_1}$. Thus, there exist $a,b\in I^D_{\xi_1}$ such that $a<t_0<b$, which implies that $g(a)>g(t_0)=0>g(b)$.

By Etemadi's strong law of large numbers (cf.\ Etemadi \cite{Ete}), for each $t\in I^D_{\xi_1}$,
 \[
  \PP\Big(\{\omega\in\Omega : \lim_{n\to\infty}g_{n,\omega}(t) = g(t) \}\Big)=1,
 \]
 and, since $\QQ$ is countable,
 \begin{align}\label{help_SLLN_deviation_1}
  \PP\Big(\{\omega\in\Omega : \lim_{n\to\infty}g_{n,\omega}(t) = g(t),\;\; \forall\; t\in I^D_{\xi_1}\cap \QQ \}\Big)=1.
 \end{align}
Let $\omega\in\Omega$ be any element such that $\lim_{n\to\infty}g_{n,\omega}(t) = g(t)$ for each $t\in I^D_{\xi_1}\cap \QQ$.
Let $\vare>0$ and $t',t''\in I^D_{\xi_1}\cap \QQ$ such that
 \[
   a\leq t_0-\vare < t' < t_0 < t'' < t_0+\vare\leq b.
 \]
Such $\vare$, $t'$ and $t''$ exist, since $t_0\in (a,b)$ and $I^D_{\xi_1}$ is an interval. Since $g$ is strictly decreasing, we have $g(t') > g(t_0)=0 > g(t'')$.
By the choices of $\omega$, $t'$ and $t''$, we have
 \[
   \lim_{n\to\infty}g_{n,\omega}(t') = g(t') \qquad \text{and}\qquad \lim_{n\to\infty}g_{n,\omega}(t'') = g(t''),
 \]
 and hence there exists $n_0:=n_0(\omega,\vare,t',t'')\in\NN$ such that $g_{n,\omega}(t')>0> g_{n,\omega}(t'')$ for all $n\geq n_0$.
Consequently, the unique zero $\mathscr{M}^D_n(\xi_1(\omega),\ldots,\xi_n(\omega))$ of $g_{n,\omega}$ is between $t'$ and $t''$ for  all $n\geq n_0$, and hence it is between $t_0-\vare$ and $t_0+\vare$ for all $n\geq n_0$. Therefore, $\mathscr{M}^D_n(\xi_1(\omega),\ldots,\xi_n(\omega))\to t_0$ as $n\to\infty$.
By \eqref{help_SLLN_deviation_1}, the statement follows.
\proofend

In the next remark we point out that the assumptions and the statement of Theorem \ref{Thm_SLLN_deviation} remain ''invariant'' for deviations that generate the same deviation means.

\begin{Rem}\label{Rem_SLLN_invariant}
Recall that, by Dar\'oczy and P\'ales \cite{DarPal82}, two deviation means on an open (nonempty) interval $I$ of $\RR$,
 generated by $D$ and $D^*$ are identical (i.e., for all $n\in\NN$, the equality $\mathscr{M}_n^{D}=\mathscr{M}_n^{D^*}$ holds on $I^n$)
 if and only if there exists a positive function $d:I\to\RR$ such that
 $D^*(x,t)=d(t)D(x,t)$, $x,t\in I$, is valid (where the equality implies the continuity of $d$ as well).
If the assumptions (i) and (ii) of Theorem \ref{Thm_SLLN_deviation} hold for $D$ and $\xi_1$, then they also hold for
 $D^*(x,t):=d(t)D(x,t)$, $x,t\in I$, and $\xi_1$ with any continuous function $d:I\to (0,\infty)$,
 since $I^{D^*}_{\xi_1} = I^{D}_{\xi_1}$, and $\EE(D^*(\xi_1,t))=0$, $t\in I$, is equivalent to $\EE(D(\xi_1,t))=0$, $t\in I$,
 yielding also that $\mathbb{M}^{D^*}(\xi_1)= \mathbb{M}^D(\xi_1)$.
Hence $\mathscr{M}_n^{D^*}(\xi_1,\ldots,\xi_n) = \mathscr{M}_n^{D}(\xi_1,\ldots,\xi_n) \as \mathbb{M}^D(\xi_1) =  \mathbb{M}^{D^*}(\xi_1)$
 as $n\to\infty$.
\proofend
\end{Rem}

\begin{Rem}
Let $I\subset\RR$ be a nondegenerate interval, let $D:I^2\to\RR$ be a deviation
 and $(\xi_n)_{n\in\NN}$ be a sequence of i.i.d.\ random variables such that $\PP(\xi_1=x_0)=1$ \ with some $x_0\in I$.
Then it can happen that $\xi_1$ does not satisfy the interior point condition with respect to $D$
 (for an example, see part (iv) of Example \ref{Ex_2}), and hence the assumption (ii) of Theorem \ref{Thm_SLLN_deviation} does not hold.
However, in this degenerate case, \eqref{help_dev_SLLN_allitas} still holds true.
Indeed, for each $n\in\NN$, $\mathscr{M}_n^{D}(\xi_1,\ldots,\xi_n) = \mathscr{M}_n^{D}(x_0,\ldots,x_0) = x_0$ almost surely,
 and $\EE(D(\xi_1,t)) = D(x_0,t)$, $t\in I$, yielding that $\mathbb{M}^D(\xi_1) = x_0$ (due to the properties of a deviation).
Consequently, $\mathscr{M}_n^{D}(\xi_1,\ldots,\xi_n) \as \mathbb{M}^D(\xi_1)$ as $n\to\infty$, as desired.
\proofend
\end{Rem}

\section{Central limit theorem for deviation means}
\label{section_CLT}

Now we present the second main result of our paper.

\begin{Thm}\label{Thm_CLT_deviation}
Let $I\subset\RR$ be a nondegenerate interval, let $D:I^2\to\RR$ be a deviation and $(\xi_n)_{n\in\NN}$ be a sequence of i.i.d.\ random variables. Assume that
\vspace{-12pt}
 \begin{enumerate}[(i)]\itemsep=-1pt
 \item $D$ is Borel measurable in the first variable;
 \item $\PP(\xi_1\in I)=1$ and $\xi_1$ satisfies the interior point condition with respect to $D$;
 \item $D$ is differentiable in the second variable at $t_0:=\mathbb{M}^D(\xi_1)$, and there exist $r>0$, $\alpha\in[0,\infty)$, a nonnegative Borel measurable function $H:I\to\RR$ and a nonnegative function $h:(0,r)\to\RR$ such that
 \begin{align}\label{D-estimate}
    \bigg|\frac{D(x, t)- D(x, t_0)}{t-t_0}- \partial_2D(x, t_0)\bigg| \leq H(x)h(|t-t_0|)
 \end{align}
 for $x\in I$, $t\ne t_0$, $t\in (t_0 - r, t_0 +r)\subset I^D_{\xi_1}$,
 and
 \[
   \limsup_{s\to0^+}\frac{h(s)}{s^\alpha}\,\,
  \begin{cases}
   =0 & \mbox{if }\, \alpha=0,\\[2mm]
   <\infty & \mbox{if }\, \alpha>0,
   \end{cases}
   \qquad\mbox{and}\qquad
   \EE(H(\xi_1)^q)<\infty,
 \]
 where $q:=\frac{2}{2+\alpha}$;
 \item $0<\EE( D(\xi_1,t_0)^2)<\infty$ and $0<\EE(-\partial_2 D(\xi_1,t_0))<\infty$.
\end{enumerate}\vspace{-12pt}
Then
 \[
    \sqrt{n}\big(\mathscr{M}_n^{D}(\xi_1,\ldots,\xi_n) - t_0\big)\distr \cN\left(0,\frac{\EE(D(\xi_1,t_0)^2 )}{( \EE(\partial_2 D(\xi_1,t_0)) )^2}\right)
             \qquad \text{as $n\to\infty$.}
 \]
 \end{Thm}

\noindent{\bf Proof.}
The facts that $D(\xi_1,t_0)$  and $\partial_2 D(\xi_1,t_0)$ are random variables are explained in the proof of Lemma \ref{Lem_aux_dev_as}, which we will use in the present proof.
Since $\PP(\xi_1\in I)=1$, we have $\PP(\xi_i\in I, i\in\NN)=1$, and hence, without loss of generality, we may assume that
 $\xi_i(\omega)\in I$ for each $\omega\in \Omega$ and $i\in\NN$.
To simplify the notations, for each $n\in\NN$, we define the random variable $\mu_n$ by
 \[
    \mu_n:=\mathscr{M}_n^D(\xi_1,\ldots,\xi_n),
 \]
 which is indeed well-defined due to Lemma \ref{Lem_Dev_mean_Borel_measurable}.
To prove the statement, we are going to show that, for all $y\in\RR$,
\begin{align}\label{help_CLT_bizo}
 \begin{split}
 &\lim_{n\to\infty}\PP\big(\sqrt{n}(\mathscr{M}_n^{D}(\xi_1,\ldots,\xi_n) - t_0)<y\big)\\
 &\qquad
    = \lim_{n\to\infty}\PP(\sqrt{n}(\mu_n - t_0)<y)
    = \Phi\bigg(-\frac{\EE(\partial_2 D(\xi_1,t_0))}{\sqrt{\EE(D(\xi_1,t_0)^2} )}y\bigg),
 \end{split}
\end{align}
where $\Phi$ stands for the distribution function of a standard normally distributed random variable.

Let $\omega\in\Omega$ and $y\in\RR$ be fixed.
By assumption (ii), $t_0$ belongs to the interior of $I^D_{\xi_1}$, thus there exists $n_0\in\NN$ (depending on $y$)
 such that $t_0+\frac{y}{\sqrt{n}}\in I_{\xi_1}^D$ for $n\geq n_0$.
Hence, using the strict monotonicity property of $D$ in the second variable, for $n\geq n_0$ we have that
\[
 \sqrt{n}(\mu_n(\omega) - t_0)<y
\]
holds if and only if $\mu_n(\omega) < t_0 + \frac{y}{\sqrt{n}}$, which is equivalent to
 \begin{align}\label{help_dev_CLT_1}
    \sum_{i=1}^n D\Big(\xi_i(\omega), t_0 + \frac{y}{\sqrt{n}} \Big)
       <   \sum_{i=1}^n D(\xi_i(\omega), \mu_n(\omega)) = 0.
 \end{align}
This inequality can be rewritten as
 \begin{align*}
  \frac{1}{\sqrt{n}}\sum_{i=1}^n D\Big(\xi_i(\omega), t_0 \Big)
  &+\frac{y}{n} \sum_{i=1}^n \partial_2D\Big(\xi_i(\omega), t_0 \Big)\\
  &+\frac{1}{\sqrt{n}}\sum_{i=1}^n \bigg(D\Big(\xi_i(\omega), t_0 + \frac{y}{\sqrt{n}} \Big)- D\Big(\xi_i(\omega), t_0 \Big)
  -\frac{y}{\sqrt{n}} \partial_2D\Big(\xi_i(\omega), t_0 \Big)\bigg)\!<0.
 \end{align*}
Here, by the first assumption in (iv) and the central limit theorem for the i.i.d.\ random variables $\big(D(\xi_i,t_0)\big)_{i\in\NN}$ (whose first moment equals zero), we obtain
 \[
    \frac{1}{\sqrt{n}}\sum_{i=1}^n D(\xi_i, t_0) \distr \cN(0, \EE(D(\xi_1,t_0)^2))
      \qquad \text{as \ $n\to\infty$.}
 \]
Further, by the second assumption in (iv) and the strong law of large numbers for the i.i.d.\ random variables $\big(\partial_2D(\xi_i, t_0)\big)_{i\in\NN}$, we get
 \begin{align}\label{help_dev_CLT_3}
   \frac{1}{n}\sum_{i=1}^n \partial_2D(\xi_i, t_0)
     \as \EE(\partial_2 D(\xi_1,t_0))\qquad \text{as \ $n\to\infty$.}
 \end{align}
Moreover, by Lemma \ref{Lem_aux_dev_as}, we have
 \begin{align}\label{help_dev_CLT_5}
 \frac{1}{\sqrt{n}}\sum_{i=1}^n \bigg(D\Big(\xi_i, t_0 + \frac{y}{\sqrt{n}} \Big)- D(\xi_i, t_0)
                                        -\frac{y}{\sqrt{n}} \partial_2D(\xi_i, t_0)\bigg)
   \as 0 \qquad \text{as \ $n\to\infty$.}
 \end{align}
Indeed, for $y=0$, it trivially holds, and for $y\ne 0$, it follows from condition (iii) and Lemma \ref{Lem_aux_dev_as} by choosing
 $t_0:=\mathbb{M}^D(\xi_1)$, $\kappa_n:=\frac{y}{\sqrt{n}}$
 and $\tau_n:=\frac{1}{\sqrt{n}}$, $n\in\NN$, since in this case
 \[
   \limsup_{n\to\infty} n^{\frac{1}{q}} \vert \tau_n\vert \vert \kappa_n\vert^{\alpha+1}
      = \limsup_{n\to\infty} \vert y\vert^{\alpha+1} n^{\frac{2+\alpha}{2}}n^{-\frac{1}{2}} n^{-\frac{\alpha+1}{2}}
      = \limsup_{n\to\infty} \vert y\vert^{\alpha+1} = \vert y\vert^{\alpha+1} <\infty.
 \]

All in all, for each $y\in\RR$, there exists $n_0\in\NN$ such that
 \[
   \PP(\sqrt{n}(\mu_n - t_0)<y)
    = \PP(\eta_n+\zeta_n<0), \qquad n\geq n_0,
 \]
 with
 \[
 \eta_n:= \frac{1}{\sqrt{n}} \sum_{i=1}^n D(\xi_i,t_0)\distr \eta \qquad \text{as $n\to\infty$},
 \]
 and
 \begin{align*}
 \zeta_n
  &:= \frac{y}{n} \sum_{i=1}^n \partial_2D\big(\xi_i, t_0 \big) +\frac{1}{\sqrt{n}}\sum_{i=1}^n \bigg(D\Big(\xi_i, t_0 + \frac{y}{\sqrt{n}} \Big)- D\big(\xi_i, t_0 \big)
       -\frac{y}{\sqrt{n}} \partial_2D\big(\xi_i, t_0 \big)\bigg) \\
  &\as \zeta \qquad \text{as $n\to\infty$},
 \end{align*}
 where $\eta$ is a random variable with distribution $\cN(0,\EE(D(\xi_1,t_0))^2)$ and $\zeta:=y\EE(\partial_2D(\xi_1, t_0))$.
By Slutsky's Lemma and by the continuity of the distribution function of $\eta+\zeta$, it follows that $\lim_{n\to\infty}\PP(\eta_n+\zeta_n<0) = \PP(\eta+\zeta<0)$.
Since
 \[
   \PP(\eta+\zeta<0)=\PP(\eta<-\zeta)  = \Phi\bigg(-\frac{\EE(\partial_2 D(\xi_1,t_0))}{\sqrt{\EE(D(\xi_1,t_0)^2} )}y\bigg),
 \]
 we conclude that \eqref{help_CLT_bizo} must hold, as desired.
\proofend

In the next remark we point out the fact that the assumptions (i), (ii) and (iv), and the statement of Theorem \ref{Thm_CLT_deviation}
 remain ''invariant'' for deviations that generate the same deviation means. We also discuss how one can generalize condition (iii) of Theorem \ref{Thm_CLT_deviation} in order to remain ''invariant'' in the same sense.

\begin{Rem}
Recall that, by Dar\'oczy and P\'ales \cite{DarPal82}, two deviation means on an open (nonempty) interval $I$ of $\RR$,
 generated by $D$ and $D^*$ are identical (i.e., for all $n\in\NN$, the equality $\mathscr{M}_n^{D}=\mathscr{M}_n^{D^*}$ holds on $I^n$)
 if and only if there exists a positive function $d:I\to\RR$ such that
 $D^*(x,t)=d(t)D(x,t)$, $x,t\in I$, is valid (where the equality also implies the continuity of $d$).
If the assumptions (i), (ii) and (iv) of Theorem \ref{Thm_CLT_deviation} hold for $D$ and $\xi_1$,
 then they also hold for $D^*(x,t):=d(t)D(x,t)$, $x,t\in I$, and $\xi_1$ with any function $d:I\to (0,\infty)$
 which is differentiable at $t_0=\mathbb{M}^{D}(\xi_1)$,
 and we also have $\EE(D^*(\xi_1,t))=0$, $t\in I$, is equivalent to $\EE(D(\xi_1,t))=0$, $t\in I$,
 yielding that $\mathbb{M}^{D^*}(\xi_1)= \mathbb{M}^D(\xi_1)$.
Indeed, for the assumptions (i) and (ii), it readily follows, and for the assumption (iv), from
 $\EE(D^*(\xi_1,t_0)^2) = d(t_0)^2\EE(D(\xi_1,t_0)^2)\in(0,\infty)$, and\
 \begin{align*}
  \EE(-\partial_2 D^*(\xi_1,t_0))
      &= \EE(-d'(t_0) D(\xi_1,t_0) - d(t_0)\partial_2 D(\xi_1,t_0))\\
      &= -d'(t_0) \EE(D(\xi_1,t_0)) + d(t_0) \EE(-\partial_2 D(\xi_1,t_0))\\
      &= d(t_0) \EE(-\partial_2 D(\xi_1,t_0)) \in(0,\infty).
 \end{align*}
Further,
 \[
   \frac{\EE(D^*(\xi_1,t_0)^2 )}{( \EE(\partial_2 D^*(\xi_1,t_0)) )^2}
      = \frac{\EE(D(\xi_1,t_0)^2 )}{( \EE(\partial_2 D(\xi_1,t_0)) )^2}.
 \]

Next we discuss condition (iii) of Theorem \ref{Thm_CLT_deviation} from the point of view ''invariance'' for deviations that generate the same deviation means.
Let us suppose the assumptions (i)-(iv) of Theorem \ref{Thm_CLT_deviation} hold for $D$ and $\xi_1$.
Then for $x\in I$, $t\ne t_0$, and $t\in(t_0-r,t_0+r)$, we have
 \begin{align*}
  \bigg|&\frac{D^*(x, t)- D^*(x, t_0)}{t-t_0}
  - \partial_2D^*(x, t_0)\bigg|\\
  &=\bigg|  \bigg(\frac{d(t)-d(t_0)}{t-t_0}-d'(t_0)\bigg)D(x,t_0)
             + d(t)\bigg(\frac{D(x,t) - D(x,t_0)}{t-t_0} - \partial_2D(x, t_0)\bigg) \\
  &\phantom{=\bigg|\,}
             + (d(t) - d(t_0))\partial_2D(x, t_0)  \bigg|\\
  &\leq \left\vert\frac{d(t)-d(t_0)}{t-t_0}-d'(t_0)\right\vert \vert D(x,t_0) \vert + d(t) H(x)h(|t-t_0|) + \vert d(t)-d(t_0)\vert \vert \partial_2D(x, t_0)\vert
 \end{align*}
 \begin{align*}
  &\leq H_1(x)h^*_1(|t-t_0|) + H_2(x)h^*_2(|t-t_0|) + H_3(x)h^*_3(|t-t_0|).
 \end{align*}
 Here
 \begin{align*}
   H_1(x):= \vert D(x,t_0) \vert , \qquad  H_2(x):=H(x),\qquad H_3(x):=\vert \partial_2D(x, t_0)\vert,\qquad  x\in I,
 \end{align*}
 and
 \[
  h_i^*(s):=\max(h_i(t_0+s),h_i(t_0-s)),\qquad s\in(0,r),\,i\in\{1,2,3\},
 \]
 where, for $t\in (t_0-r,t_0+r)$ with $t\neq t_0$,
  \[
  h_1(t):=\left\vert\frac{d(t)-d(t_0)}{t-t_0}-d'(t_0)\right\vert,\qquad
  h_2(t):=d(t)h(|t-t_0|),\qquad
  h_3(t):=\vert d(t)-d(t_0)\vert.
 \]
Indeed, for $t\in (t_0-r,t_0+r)$ with $t\neq t_0$, we have
  \[
  h_i(t)\leq \max(h_i(t),h_i(2t_0-t))=\max(h_i(t_0+|t-t_0|),h_i(t_0-|t-t_0|))=h_i^*(|t-t_0|).
 \]
Here, by conditions (iii) and (iv) of Theorem \ref{Thm_CLT_deviation}, for $i\in\{1,2,3\}$, we have that $\EE(H_i^{q_i}(\xi_1))<\infty$, where $(q_1,q_2,q_3)=(1,q,1)$. Furthermore, by the differentiability and continuity of $d$ at $t_0$, respectively, we have that
\[
 \limsup_{s\to0^+}h_1^*(s)=\limsup_{s\to0^+}h_3^*(s)=0,
\]
and
\[
 \limsup_{s\to0^+}\frac{h_2^*(s)}{s^\alpha}
 =d(t_0)\limsup_{s\to0^+}\frac{h(s)}{s^\alpha}.
\]

The above argument shows that in general the condition (iii) of Theorem \ref{Thm_CLT_deviation} is not invariant for deviations that generate the same deviation means. Applying the method of the proof of Theorem~\ref{Thm_CLT_deviation}, one can verify that the conclusion of this theorem remains valid if its condition (iii) is replaced by the following more general one:
 \begin{enumerate}[(i)]\itemsep=-4pt
 \item[(iii$^*$)] $D$ is differentiable in the second variable at $t_0:=\mathbb{M}^D(\xi_1)$, and there exist $r>0$, $k\in\NN$, $\alpha_1,\dots,\alpha_k\in[0,\infty)$, nonnegative Borel measurable functions $H_1,\dots,H_k:I\to\RR$ and nonnegative functions $h_1,\dots,h_k:(0,r)\to\RR$ such that
 \begin{align*}
    \bigg|\frac{D(x, t)- D(x, t_0)}{t-t_0}- \partial_2D(x, t_0)\bigg| \leq \sum_{i=1}^kH_i(x)h_i(|t-t_0|)
 \end{align*}
 for $x\in I$, $t\ne t_0$, $t\in (t_0 - r, t_0 +r)\subset I^D_{\xi_1}$,
 and for $i\in\{1,\dots,k\}$,
 \[
   \limsup_{s\to0^+}\frac{h_i(s)}{s^{\alpha_i}}\,\,
  \begin{cases}
   =0 & \mbox{if }\, \alpha_i=0,\\[2mm]
   <\infty & \mbox{if }\, \alpha_i>0,
   \end{cases}
   \qquad\mbox{and}\qquad
   \EE(H(\xi_1)^{q_i})<\infty,
 \]
 where $q_i:=\frac{2}{2+\alpha_i}$.
\end{enumerate}\vspace{-12pt}
\proofend
\end{Rem}

\begin{Rem}
Note that the degenerate case when $\PP(\xi_1=x_0)=1$ \ with some $x_0\in I$
 is excluded from Theorem \ref{Thm_CLT_deviation}, since in this case
  $\EE(D(\xi_1,t_0)^2) = \EE(D(x_0,t_0)^2) = \EE(0^2) =0$ and
 it can happen that $\xi_1$ does not satisfy the interior point condition
 with respect to $D$ (for an example, see part (iv) of Example \ref{Ex_2}),
 and hence the assumption (iv) (and maybe (ii)) of Theorem \ref{Thm_CLT_deviation} does not hold.
However, in this degenerate case, for each $n\in\NN$ we have
 \[
    \sqrt{n}\big(\mathscr{M}_n^{D}(\xi_1,\ldots,\xi_n) - t_0\big)
      = \sqrt{n}\big(\mathscr{M}_n^{D}(x_0,\ldots,x_0) - x_0\big)
      = \sqrt{n}(x_0 - x_0)
      = 0
 \]
 almost surely, and thus one can say that the statement of Theorem \ref{Thm_CLT_deviation} still holds true
 provided that $\cN(0,0)$ is understood as the Dirac measure on $\RR$ concentrated at the point $0$.
\proofend
\end{Rem}

For establishing a parameter $\alpha$ which satisfies condition (iii) of Theorem~\ref{Thm_CLT_deviation}, the following assertion could be useful.

\begin{claim} Let $r>0$ and $h:(0,r)\to(0,\infty)$. If there exists $\alpha\in[0,\infty)$ such that
 \begin{align}\label{help_alpha_limsup0}
 \limsup_{s\to0^+}\frac{h(s)}{s^\alpha}<\infty,
 \end{align}
then
 \begin{align}\label{help_alpha_limsup1}
   \alpha\leq \liminf_{s\to0^+}\frac{\ln(h(s))}{\ln(s)}.
 \end{align}
Conversely, if $\alpha$ satisfies \eqref{help_alpha_limsup1} with a strict inequality, then \eqref{help_alpha_limsup0} holds.
\end{claim}

\noindent{\bf Proof.}
Assume first that $A:=\limsup_{s\to0^+}\frac{h(s)}{s^\alpha}<\infty$. Then there exists $s_1\in(0,\min(1,r))$ such that
 \[
    0<\frac{h(s)}{s^\alpha} < A+1, \qquad 0<s<s_1,
 \]
 and hence
 \[
    \ln(h(s)) < \ln(A+1) + \alpha\ln(s), \qquad 0<s<s_1,
 \]
 that is,
 \[
    \frac{\ln(h(s))}{\ln(s)} > \alpha + \frac{\ln(A+\vare)}{\ln(s)}, \qquad 0<s<s_1,
 \]
 which, upon taking the limit $s\to0^+$, yields \eqref{help_alpha_limsup1}.

 Conversely, suppose that $\alpha$ satisfies \eqref{help_alpha_limsup1} with a strict inequality.  Then there exists $s_2\in(0,\min(1,r))$ such that
 \[
    \frac{\ln(h(s))}{\ln(s)} > \alpha, \qquad 0<s<s_2,
 \]
 which implies that
 \[
    0<\frac{h(s)}{s^\alpha} < 1, \qquad 0<s<s_2.
 \]
 By taking the limit $s\to0^+$, this inequality yields that \eqref{help_alpha_limsup0} is valid.
\proofend

In the next remark we provide a sufficient condition under which the assumption (ii) of Theorem \ref{Thm_CLT_deviation} holds.
In particular, the given argument also shows that in case of $\alpha=0$, the second part of assumption (ii) of Theorem \ref{Thm_CLT_deviation},
 i.e., the assumption that $\xi_1$ satisfies the interior point condition with respect to $D$, holds provided that
 $\PP(\xi_1\in I)=1$ and the assumptions (i), (iii) and (iv) of Theorem \ref{Thm_CLT_deviation} hold with
 a point $t_0:=\mathbb{M}^D(\xi_1)$ satisfying $\EE(D(\xi_1,t_0))=0$  and being an interior point of $I$.

\begin{Rem}
Let $I\subset\RR$ be a nondegenerate interval, let $D:I^2\to\RR$ be a deviation and $(\xi_n)_{n\in\NN}$ be a sequence of i.i.d.\ random variables.
Assume that $\PP(\xi_1\in I)=1$, there exists a point $t_0:=\mathbb{M}^D(\xi_1)$ such that $\EE(D(\xi_1,t_0))=0$ and
 $t_0$ is an interior point of $I$.
Further, assume that the assumptions (i) and (iv) of Theorem \ref{Thm_CLT_deviation} hold,
  $D$ is differentiable in the second variable at $t_0$ and there exist $r>0$, a nonnegative Borel measurable function
  $H:I\to\RR$ and nonnegative function $h:(0,r)\to\RR$ such that
 \begin{align}\label{D-estimate_2}
    \bigg|\frac{D(x, t)- D(x, t_0)}{t-t_0}- \partial_2D(x, t_0)\bigg| \leq H(x)h(|t-t_0|)
 \end{align}
 for $x\in I$, $t\ne t_0$, $t\in (t_0 - r, t_0 +r)\subset I$ such that $\EE(H(\xi_1))<\infty$.
Then $\xi_1$ satisfies the interior point condition with respect to $D$ yielding that the assumption (ii) of Theorem \ref{Thm_CLT_deviation} holds.
Indeed, by \eqref{D-estimate_2},
 \[
   \big|D(x, t)\big| \leq  \big|D(x, t_0)\big| + |t-t_0| \big( |\partial_2D(x, t_0)| + H(x)h(|t-t_0|) \big)
 \]
 for $x\in I$, $t\ne t_0$ and $t\in (t_0 - r, t_0 +r)$.
Consequently, using that $\EE(\big|D(\xi_1, t_0)\big|)<\infty$, $\EE(|\partial_2D(\xi_1, t_0)|)<\infty$ and $\EE(H(\xi_1))<\infty$
 (due to the assumptions), we have $\EE( \big|D(\xi_1, t)\big|)<\infty$ for $t\in (t_0 - r, t_0 +r)$.
Hence $(t_0 - r, t_0 +r)\subset I^D_{\xi_1}$, and thus $t_0$ is an interior point of $I^D_{\xi_1}$.
\proofend
\end{Rem}

In the next remark we provide some sufficient conditions under which the inequality \eqref{D-estimate} in assumption (iii) of Theorem \ref{Thm_CLT_deviation} holds.

\begin{Rem}\label{Rem_CLT_a(iii)_elegseges}
(i).\
Suppose that $I\subset \RR$ is a nondegenerate interval, $t_0\in I$, $D:I^2\to\RR$ is a deviation which is differentiable in the second variable and there exist $r>0$, a nonnegative function $H:I\to\RR$ and a nondecreasing nonnegative function $h:[0,r)\to\RR$ such that
 \begin{align}\label{help_assump_iv_ver1}
    |\partial_2 D(x,t)-\partial_2 D(x,t_0)| \leq H(x)h(|t-t_0|),\qquad x\in I, \;\; t\in (t_0 - r, t_0 +r)\subset I.
 \end{align}
Then the inequality \eqref{D-estimate} in assumption (iii) of Theorem \ref{Thm_CLT_deviation} is satisfied for $x\in I$ and $t\ne t_0$, $t\in (t_0 - r, t_0 +r)$.
Indeed, for $x\in I$ and $t\ne t_0$, $t\in (t_0 - r, t_0 +r)$, by the Lagrange mean value theorem, there exists a value $t^*$ between $t$ and $t_0$ (depending on $x$) such that
 \begin{align*}
  \bigg|\frac{D(x, t)- D(x, t_0)}{t-t_0}- \partial_2D(x, t_0)\bigg|
  & = \vert \partial_2D(x, t^*) - \partial_2D(x, t_0) \vert \\
  & \leq H(x) h(\vert t^*-t_0\vert)
   \leq  H(x) h(\vert t-t_0\vert),
 \end{align*}
 where we used \eqref{help_assump_iv_ver1} and the monotonicity of $h$.

\noindent (ii).\ By Rademacher's theorem, if $h(t)=t$, $t\in[0,r)$ in \eqref{help_assump_iv_ver1}, then for each $x\in I$, the function $t\ni (t_0-r,t_0+r)\mapsto \partial_2 D(x,t)$ is (globally) Lipschitz continuous,
 and hence it is differentiable at Lebesgue almost every $t\in (t_0-r,t_0+r)$.
Consequently, in this case for each $x\in I$, we have that $\partial_2^2 D(x,t)$ exists at Lebesgue almost every $t\in (t_0-r,t_0+r)$.
In part (iii) of the remark, we formulate a condition on $\partial_2^2 D$ under which \eqref{help_assump_iv_ver1} holds with $h(t):=t$, $t\in[0,r)$.

\noindent (iii).\ Assume that $I\subset \RR$ is a nondegenerate interval, $t_0\in I$, $D:I^2\to\RR$ is a deviation which is twice differentiable
 with respect to the second variable, and there exist $r>0$ and a nonnegative function $H:I\to\RR$ such that
 \begin{align}\label{help_assump_iv_ver2}
    \vert \partial_2^2 D(x,t)\vert \leq H(x),\qquad x\in I, \;\; t\in(t_0 - r, t_0 +r) \subset I.
 \end{align}
Then \eqref{help_assump_iv_ver1} (and hence the inequality \eqref{D-estimate} in assumption (iii) of Theorem \ref{Thm_CLT_deviation} due to part (i) of the remark)
 holds with $h(t):=t$, $t\in[0,r)$. For $t=t_0$, the statement is trivial.
If $t\in(t_0-r,t_0+r)$ and $t\ne t_0$, then, by the Lagrange's mean value theorem, there exists a value $t^*$ (depending on $x$)
 between $t$ and $t_0$ such that
 \begin{align}\label{help_assump_iv_ver3}
   |\partial_2 D(x,t)-\partial_2 D(x,t_0)|
     = |\partial_2^2 D(x,t^*)||t-t_0|
     \leq H(x) |t-t_0|,\qquad x\in I.
 \end{align}
This proves the assertion.

\noindent (iv).\
In the cases when the inequality \eqref{D-estimate} in assumption (iii) of Theorem \ref{Thm_CLT_deviation} holds with $h(t):=t$, $t\in(0,r)$, the second part of this assumption is satisfied with $\alpha=1$, thus, $q=\frac{2}{3}$ and only $\EE(H(\xi_1)^{\frac{2}{3}})<\infty$ is required.
\proofend
\end{Rem}

In the next remark we compare the assumptions of Theorems \ref{Thm_SLLN_deviation} and \ref{Thm_CLT_deviation} with the assumptions of \emph{some} existing limit theorems for maximum likelihood estimators and M-estimators.
There are plenty sets of sufficient assumptions under which limit theorems for M-estimators hold, and in our forthcoming remark we point out that our set of sufficient conditions for deviation means are different from the ones available for M-estimators in the literature, and we think that they are more simpler and more checkable.

\begin{Rem}\label{Rem_feltet_osszehas}
\noindent (i).\
In limit theorems on asymptotic normality of maximum likelihood estimators of some parameters based on i.i.d.\ observations
 somewhat similar assumptions as \eqref{help_assump_iv_ver1} and \eqref{help_assump_iv_ver2} appear.
Inequality \eqref{help_assump_iv_ver2} with $\EE(H(\xi_1))<\infty$ can be considered as a counterpart of the conditions (3.15) and (3.16) in Lehmann and Casella \cite[Theorem 6.3.10]{LehCas98}.
Inequality \eqref{help_assump_iv_ver3} with $\EE(H(\xi_1)^2)<\infty$ can be considered as a counterpart of the condition in Theorem 7.12 in van der Vaart \cite{Vaa}.
Finally, note that \eqref{help_assump_iv_ver1} with  $h(t):=t^\beta$, $t\in[0,r)$, where $\beta\in(0,1]$, can be considered as a counterpart of the uniform Lipschitz type condition (c) of Theorem 2 in Rao \cite[Chapter 3, Section 4]{Rao}.

\indent (ii).\ It was first Huber \cite[page 74]{Hub64} who introduced the concept of an M-estimator, namely, using the notations introduced in Remark \ref{Rem_M_estimate}, he considered the case when $X:=\RR$, $\Theta=\RR$,  and the function $\varrho$ depends only on $x-\vartheta$, i.e., $\varrho(x,\vartheta):=f(x-\vartheta)$, $x\in \RR$, $\vartheta\in\Theta$, with some given nonconstant function $f:\RR\to\RR$. It is in fact a special M-estimator. The asymptotic behaviour of this M-estimator based on i.i.d.\ observations has been investigated according to the two cases, namely, when $f$ is convex and nonconvex. Provided that $f$ is continuous, convex and it has limit $\infty$ at $\pm\infty$, the result
            of Huber \cite[Lemma 1]{Hub64} states that the set of the solutions of the minimization problem \eqref{help_M_est_min_problem}
            is nonempty, convex and compact, and there is a unique minimizer provided that $f$ is strictly convex.
            Provided that $f$ is differentiable, Lemma 2 in Huber \cite{Hub64} can be considered
            as a counterpart of our Lemma \ref{Lem_Dev_mean_random_1}.
            Provided that $f$ is differentiable, Lemma 3 in Huber \cite{Hub64} provides strong consistency
            of the M-estimator in question under the condition that there exists a $c\in\RR$ such that
            $\EE(f'(\xi_1-t))<0$ for $t>c$ and $\EE(f'(\xi_1-t))>0$ for $t<c$.
            This later condition might be considered as a counterpart of our interior point condition assumed in Theorem \ref{Thm_SLLN_deviation},
            however, note that we do not assume any differentiability for $D$ in Theorem \ref{Thm_SLLN_deviation}.
            Provided that $f$ is differentiable, Lemma 4 in Huber \cite{Hub64} provides asymptotic normality of the M-estimator in question
            under some sufficient conditions, which seem to be stronger than the ones in our Theorem \ref{Thm_CLT_deviation}.
            For example, it is assumed that the function $\RR\ni\vartheta\mapsto \EE( (f'(\xi_1-\vartheta))^2 )$ is finite and continuous
            at the point $c$ given before, and note that we do not have such a condition on $D$ and $\xi_1$ in Theorem \ref{Thm_CLT_deviation}.
            Nonetheless, the first parts of the proofs of Lemma 4 in Huber \cite{Hub64} and our Theorem \ref{Thm_CLT_deviation}
            are somewhat similar.
            If $f$ is not convex, than Huber \cite[page 78]{Hub64} noted that, in general, strong consistency of the M-estimator in question
            does not hold.
            Further, in case of a nonconvex $f$, provided that the M-estimator in question converges in probability towards 0 as the sample size tends to 0,
            Huber \cite[Lemma 5]{Hub64} gives a set of sufficient conditions under which asymptotic normality of
            the M-estimator in question holds.
            Among others, these sufficient conditions include the uniform continuity of the second derivative of $f$,
            and note that we do not assume anything similar in our Theorem \ref{Thm_CLT_deviation}.

\indent (iii).\ Huber \cite{Hub67} introduced a generalization of M-estimators compared to the one recalled in Remark \ref{Rem_M_estimate}.
          Namely, keeping the notations introduced in Remark \ref{Rem_M_estimate},
          a sequence $\widehat\vartheta_n:=\widehat\vartheta_n(\xi_1,\ldots,\xi_n)$, $n\in \NN$,
          can be called a(n approximative) M-estimator if
         \begin{align*}
         \frac{1}{n} \sum_{i=1}^n \varrho(\xi_i,\widehat\vartheta_n)
            - \argmin_{\vartheta\in\Theta}\frac{1}{n} \sum_{i=1}^n \varrho(\xi_i,\vartheta)
            \to 0 \qquad \text{as \ $n\to\infty$ \ almost surely.}
        \end{align*}
        In case of a locally compact parameter set $\Theta$, Theorem 1 in Huber \cite{Hub67} contains some
        sufficient conditions under which strong consistency of the approximative M-estimator in question holds.
        Huber \cite[page 224]{Hub67} also introduced a kind of approximative $\psi$-estimator, namely,
        given a function $\psi:X\times\Theta\to\RR$, a sequence
        $\widehat\vartheta_n:=\widehat\vartheta_n(\xi_1,\ldots,\xi_n)$, $n\in \NN$, can be called a(n approximative) $\psi$-estimator if
        \[
         \frac{1}{n}\sum_{i=1}^n \psi(\xi_i,\widehat\vartheta_n)\to 0 \qquad \text{as \ $n\to\infty$ \ almost surely.}
        \]
        The connection between approximative M-estimators and approximative $\psi$-estimators might be studied via Ekeland's variational principle \cite{Eke74}.
        In case of a locally compact parameter set $\Theta$, Theorem 2 in Huber \cite{Hub67} contains some sufficient conditions
        under which strong consistency of the approximative $\psi$-estimator in question holds.
        These sufficient conditions seem to be somewhat different that we assume in our Theorem \ref{Thm_SLLN_deviation}.
        For example, Huber \cite[Theorem 2]{Hub67} assumes that for each $\vartheta\in\Theta$, we have
        $\EE(\sup_{\vartheta'\in U} \vert \psi(\xi_1,\vartheta') - \psi(\xi_1,\vartheta) \vert)  \to 0$ as the ball $U$
        around $\vartheta$ shrinks to $\vartheta$, and note that we do not have such a condition in our Theorem \ref{Thm_SLLN_deviation}.
        In case of an open parameter set $\Theta$, Huber \cite[Corollary on page 231]{Hub67} gives a set of sufficient conditions under which
        asymptotic normality holds for the approximative $\psi$-estimator in question.
        These sufficient conditions seem to be quite different from the assumptions in our Theorem \ref{Thm_CLT_deviation}.
        For example, we do not have any assumption corresponding to the assumption (N-3)/(iii) or the differentiability
        of the function $\Theta\ni\vartheta\mapsto \EE(\psi(\xi_1,\vartheta))$ in Huber \cite[Corollary on page 231]{Hub67}.

\indent (iv).\ van der Vaart \cite[Section 5]{Vaa} also gives a detailed exposition of M-estimators and $\psi$-estimators,
         and several sets of sufficient conditions under which weak consistency and asymptotic normality hold.
        Here we only mention that Theorem 5.23 in van der Vaart \cite{Vaa} provides asymptotic normality
        of M-estimators under some sufficient conditions which are stronger than the assumptions of
        our Theorem \ref{Thm_CLT_deviation} for deviations means.
        More precisely, using the notations introduced in Remark \ref{Rem_M_estimate},
        in case of $X=\RR$, an open parameter set $\Theta$, van der Vaart \cite[Theorem 5.23]{Vaa}
        assumes for example that for each $x\in\RR$, the function $\Theta\ni\vartheta \mapsto \varrho(x,\vartheta)$ is Lipschitz continuous
        in some neighbourhood of a maximizer $\vartheta_0$ of the function $\Theta\ni\vartheta \mapsto \EE(\varrho(\xi_1,\vartheta))$
        and the function $\Theta\ni\vartheta \mapsto \EE(\varrho(\xi_1,\vartheta))$ admits a second-order Taylor expansion at $\vartheta_0$.
        Note that we do not have any corresponding assumptions in Theorem \ref{Thm_CLT_deviation}.

\indent (v).\  He and Wang \cite[Theorem 2.1 and Corollary 2.1]{HeWan} proved
              a law of the iterated logarithm for $\psi$-estimators based on i.i.d.\ observations.
              Their sufficient conditions are stronger than the assumptions of our Theorem \ref{Thm_LIL_deviation} for
              deviation means, see, e.g., the assumption (M3)/(iii) in Theorem 2.1 in \cite{HeWan} or the Lipschitz condition
              in Corollary 2.1 in \cite{HeWan}.

\indent (vi).\  Rubin and Rukhin \cite[Theorem 2]{RubRuk} proved a large deviation theorem for approximative $\psi$-estimators
                based on i.i.d.\ observations.
                Arcones \cite{Arc} established very general large deviation results for M-estimators based on i.i.d.\ observations.
\proofend
\end{Rem}

\section{Law of the iterated logarithm for deviation means}
\label{section_LIL}

Next we present the third main result of our paper.

\begin{Thm}\label{Thm_LIL_deviation}
Let $I\subset\RR$ be a nondegenerate interval, let $D:I^2\to\RR$ be a deviation and $(\xi_n)_{n\in\NN}$ be a sequence of i.i.d.\ random variables.
Let us suppose that the assumptions (i), (ii) and (iv) of Theorem \ref{Thm_CLT_deviation} on $D$, $\xi_1$ and $t_0:=\mathbb{M}^D(\xi_1)$ hold. Further, let us suppose that
 \begin{itemize}
 \item[(iii')] $D$ is differentiable in the second variable at $t_0$ and there exist $r>0$, $\alpha\in[0,\infty)$, a nonnegative Borel measurable function $H:I\to\RR$ and a nonnegative function $h:(0,r)\to\RR$ such that
 \begin{align}\label{D-estimate_LIL}
    \bigg|\frac{D(x, t)- D(x, t_0)}{t-t_0}- \partial_2D(x, t_0)\bigg| \leq H(x)h(|t-t_0|)
 \end{align}
 for $x\in I$, $t\ne t_0$, $t\in (t_0 - r, t_0 +r)\subset I^D_{\xi_1}$,
 and
 \[
   \limsup_{s\to0^+}\frac{h(s)}{s^\alpha}\,\,
  \begin{cases}
   =0 & \mbox{if }\, \alpha=0,\\[2mm]
   <\infty & \mbox{if }\, \alpha>0,
   \end{cases}
   \qquad\mbox{and}\qquad
   \EE(H(\xi_1)^q)<\infty,
 \]
 where $q:=1$ if $\alpha=0$, and $q\in\big(\frac{2}{2+\alpha},1\big)$ if $\alpha>0$.
 \end{itemize}
Then we have
 \begin{align}\label{help_LIL_1}
    \limsup_{n\to\infty} \frac{\mathscr{M}_n^{D}(\xi_1,\ldots,\xi_n) - t_0}{\sqrt{2n^{-1}\ln(\ln(n))}}
       \ase  \frac{\sqrt{\EE(D(\xi_1,t_0)^2)}}{ \EE(-\partial_2 D(\xi_1,t_0))}
 \end{align}
 and
 \begin{align}\label{help_LIL_2}
    \liminf_{n\to\infty} \frac{\mathscr{M}_n^{D}(\xi_1,\ldots,\xi_n) - t_0}{\sqrt{2n^{-1}\ln(\ln(n))}}
       \ase  -\frac{\sqrt{\EE(D(\xi_1,t_0)^2)}}{ \EE(-\partial_2 D(\xi_1,t_0))}.
 \end{align}
 \end{Thm}

Note that condition (iii) of Theorem \ref{Thm_CLT_deviation} and condition (iii') of Theorem \ref{Thm_LIL_deviation} are the same if $\alpha=0$. However, they are different if $\alpha>0$, since then $q=\frac{2}{2+\alpha}$ in condition (iii), while
$q\in\big(\frac{2}{2+\alpha},1\big)$ in condition (iii').

\noindent{\bf Proof of  Theorem \ref{Thm_LIL_deviation}.}
Since $\PP(\xi_1\in I)=1$, we have $\PP(\xi_i\in I, i\in\NN)=1$, and hence, without loss of generality, we may assume that
 $\xi_i(\omega)\in I$ for each $\omega\in \Omega$ and $i\in\NN$.
To simplify the notations, let
 \[
  C:=\frac{\sqrt{\EE(D(\xi_1,t_0)^2)}}{ \EE(-\partial_2 D(\xi_1,t_0))},
  \qquad \text{and} \qquad a_n:= \sqrt{2n^{-1}\ln(\ln(n))}, \qquad n\in\NN\setminus\{1,2\}.
 \]
Recall also the notation from the proof of Theorem \ref{Thm_CLT_deviation}:
 \[
    \mu_n:=\mathscr{M}_n^D(\xi_1,\ldots,\xi_n),\qquad n\in\NN.
 \]
We introduce the following sequences of random variables:
 \begin{align*}
   &R_n:= \frac{1}{\sqrt{2n\ln(\ln(n))}} \sum_{i=1}^n D(\xi_i, t_0),\qquad n\in\NN\setminus\{1,2\},\\[1mm]
   &S_n:=\frac{1}{n} \sum_{i=1}^n \partial_2 D(\xi_i, t_0),\qquad n\in\NN,
 \end{align*}
 and
 \begin{align*}
   T_n^{(\gamma)}:=\frac{1}{\sqrt{2n\ln(\ln(n))}} \sum_{i=1}^n \Big( D(\xi_i, t_0+\gamma a_n) - D(\xi_i,t_0)- \gamma a_n \partial_2 D(\xi_i, t_0) \Big)
 \end{align*}
 for $n\in\NN\setminus\{1,2\}$, $\gamma\in\RR$ satisfying $\vert \gamma a_n\vert<r$, where $r$ appears in assumption (iii').
Note that, by assumption (i), $R_n$, $n\in\NN\setminus\{1,2\}$, are indeed random variables, and, by the proof of Lemma \ref{Lem_aux_dev_as}, $S_n$, $n\in\NN$, and $T_n^{(\gamma)}$ for $n\in\NN\setminus\{1,2\}$, $\gamma\in\RR$ satisfying $\vert \gamma a_n\vert<r$ are also random variables.

Since $\EE(D(\xi_1, t_0))=0$, the first part of assumption (iv) and the law of the iterated logarithm for the sequence of i.i.d.\ random variables
 $(D(\xi_n,t_0))_{n\in\NN}$ yield that
 \begin{align}\label{help_LIL_5}
   \limsup_{n\to\infty} R_n
    \ase \sqrt{\EE(D(\xi_1,t_0)^2)}.
 \end{align}
Further, by the second part of assumption (iv) and Kolmogorov's strong law of large numbers for the sequence of
 i.i.d.\ random variables $(\partial_2D(\xi_n,t_0))_{n\in\NN}$, we have
 \begin{align}\label{help_LIL_6}
     S_n \as \EE(\partial_2 D(\xi_1,t_0)) \qquad \text{as \ $n\to\infty$.}
 \end{align}
Moreover, by assumption (iii') and Lemma \ref{Lem_aux_dev_as}, for each $\gamma\in\RR$, we check that
 \begin{align}\label{help_LIL_12}
     T_n^{(\gamma)} \as 0 \qquad \text{as \ $n\to\infty$.}
 \end{align}
In case of $\gamma=0$, \eqref{help_LIL_12} is trivial.
In case of $\gamma\ne 0$, since $a_n\to 0$ as $n\to\infty$, we have $\vert \gamma a_n\vert<r$ for sufficiently
 large $n\in\NN$, and hence $T_n^{(\gamma)}$ is well-defined for sufficiently large $n\in\NN$.
In case of $\gamma\ne 0$, let us apply Lemma \ref{Lem_aux_dev_as} with the following choices: $t_0:=\mathbb{M}^D(\xi_1)$,
 $\kappa_n:=\gamma a_n$ and $\tau_n:=\frac{1}{\sqrt{2n\ln(\ln(n))}}$, $n\in\NN\setminus\{1,2\}$, for which
 condition (c) of Lemma \ref{Lem_aux_dev_as} holds, since
 \[
   \limsup_{n\to\infty} n^{\frac{1}{q}} \vert \tau_n\vert \vert \kappa_n\vert^{\alpha+1}
      = \begin{cases}
          \limsup\limits_{n\to\infty} \vert \gamma\vert = \vert \gamma\vert <\infty & \text{if $\alpha=0$,}\\[2mm]
          \limsup\limits_{n\to\infty} 2^{\frac{\alpha}{2}} \vert \gamma\vert^{\alpha+1} n^{\frac{1}{q} - \frac{\alpha}{2} - 1} (\ln(\ln(n)))^{\frac{\alpha}{2}}
                                    =0<\infty & \text{if $\alpha>0$,}\\
        \end{cases}
 \]
 where we used that $q>\frac{2}{2+\alpha}$ yields that $\frac{1}{q} - \frac{\alpha}{2} - 1<0$ in case of $\alpha>0$.

Using that the intersection of countably many events having probability one is an event having probability one,
 by \eqref{help_LIL_5}, \eqref{help_LIL_6} and \eqref{help_LIL_12}, for each $\gamma\in\RR$, we have $\PP(\Omega^{(\gamma)})=1$,
 where
 \begin{align*}
   \Omega^{(\gamma)}:=\Big\{
                        \omega\in\Omega :\, & \limsup_{n\to\infty} R_n (\omega) = \sqrt{\EE(D(\xi_1,t_0)^2)};\\
                                            &  \lim_{n\to\infty} S_n(\omega) = \EE(\partial_2 D(\xi_1,t_0));\;
                                               \lim_{n\to\infty} T_n^{(\gamma)}(\omega) = 0
                       \Big\}.
 \end{align*}

First, we prove \eqref{help_LIL_1}.
Let $\vare>0$ be fixed.
By assumption (ii), $t_0$ belongs to the interior of $I^D_{\xi_1}$, thus, using that $\lim_{n\to\infty} a_n=0$,  we have
there exists $n_0\in\NN\setminus\{1,2\}$ such that $t_0+(C+\vare)a_n\in I_{\xi_1}^D$ for $n\geq n_0$.
Using the strict monotonicity property of $D$ in the second variable,
for each $\omega\in\Omega$ and $n\geq n_0$, we have that the following inequalities
 \begin{align}
   \frac{\mu_n(\omega)-t_0}{a_n} &< C+\vare, \label{help_LIL_2.5}\\
   \mu_n(\omega) &< t_0 + (C+\vare)a_n,\nonumber\\
  \sum_{i=1}^n D(\xi_i(\omega), t_0+(C+\vare)a_n)
        &<  \sum_{i=1}^n D(\xi_i(\omega), \mu_n(\omega))=0,\nonumber\\
   \label{help_LIL_3}
     R_n(\omega)&< -(C+\vare)S_n(\omega) -  T_n^{(C+\vare)}(\omega)
 \end{align}
are equivalent to each other.

Let $\omega\in\Omega^{(C+\vare)}$ be fixed arbitrarily. Then we have
 \begin{align}\label{help_LIL_7}
  \begin{split}
    &-(C+\vare)S_n(\omega) - T_n^{(C+\vare)}(\omega)\\
    &\qquad \to -(C+\vare) \EE(\partial_2 D(\xi_1,t_0))
            = \sqrt{\EE(D(\xi_1,t_0)^2)} + \vare \EE(-\partial_2 D(\xi_1,t_0))
  \end{split}
 \end{align}
as $n\to\infty$. Thus, in view of the positivity of $\EE(-\partial_2 D(\xi_1,t_0))$ (according to the second part of assumption (iv)), there exists $n_\vare^{*}(\omega)\in\NN$ with $n_\vare^{*}(\omega)\geq n_0$ such that, for all $n\geq n_\vare^{*}(\omega)$, we have
 \begin{align*}
    -(C+\vare)S_n(\omega) -  T_n^{(C+\vare)}(\omega)
      &> \sqrt{\EE(D(\xi_1,t_0)^2)} + \frac{\vare}{2} \EE(-\partial_2 D(\xi_1,t_0)).
 \end{align*}
On other hand, there exists $n_\vare^{**}(\omega)\in\NN$ with $n_\vare^{**}(\omega)\geq n_0$ such that, for all $n\geq n_\vare^{**}(\omega)$, we also have
 \begin{align*}
  R_n(\omega) < \sqrt{\EE(D(\xi_1,t_0)^2)}
       +\frac{\vare}{2} \EE(-\partial_2 D(\xi_1,t_0)).
 \end{align*}
Consequently, the inequality \eqref{help_LIL_3} holds for all $n\geq n_\vare(\omega):=\max(n_\vare^{*}(\omega),n_\vare^{**}(\omega))$, which implies the validity of \eqref{help_LIL_2.5} for all $n\geq n_\vare(\omega)$. This shows that, for each $\vare>0$ and $\omega\in\Omega^{(C+\vare)}$, we have
 \[
    \limsup_{n\to\infty} \frac{\mu_n(\omega)-t_0}{a_n} \leq  C+\vare.
 \]
Then, by choosing $\vare:=\frac{1}{m}$, $m\in\NN$,  we get $\limsup_{n\to\infty} \frac{\mu_n(\omega)-t_0}{a_n} \leq  C$
 for each $\omega\in \bigcap_{m\in\NN}\Omega^{(C+\frac{1}{m})}$, where $\PP(\bigcap_{m\in\NN}\Omega^{(C+\frac{1}{m})})=1$,
 yielding that $\limsup_{n\to\infty} \frac{\mu_n-t_0}{a_n} \leq  C$ almost surely.

Let $\vare\in(0,C)$ be fixed.
By assumption (ii), $t_0$ belongs to the interior of $I^D_{\xi_1}$, thus, using that $\lim_{n\to\infty} a_n=0$,  we have
 there exists $n_1\in\NN\setminus\{1,2\}$ such that $t_0+(C-\vare)a_n\in I_{\xi_1}^D$ for $n\geq n_1$.
Using the strict monotonicity property of $D$ in the second variable and using a completely analogous argument as above, for each $\omega\in\Omega$ and $n\geq n_1$,
 we can see that the following two inequalities
 \begin{align}
   \frac{\mu_n(\omega)-t_0}{a_n} &> C-\vare, \label{help_LIL_7.5}\\
   \label{help_LIL_8}
     R_n(\omega)&> -(C-\vare)S_n(\omega) -  T_n^{(C-\vare)}(\omega)
 \end{align}
are equivalent to each other.


Let $\omega\in \Omega^{(C-\vare)}$ be fixed arbitrarily. Then we have
 \begin{align}\label{help_LIL_10}
  \begin{split}
   &-(C-\vare)S_n(\omega) - T_n^{(C-\vare)}(\omega)\\
   &\qquad \to -(C-\vare) \EE(\partial_2 D(\xi_1,t_0))
                         = \sqrt{\EE(D(\xi_1,t_0)^2)} - \vare \EE(-\partial_2 D(\xi_1,t_0))
  \end{split}
 \end{align}
 as $n\to\infty$.  Thus, in view of the positivity of $\EE(-\partial_2 D(\xi_1,t_0))$, there exists $m_\vare^{*}(\omega)\in\NN$ with $m_\vare^{*}(\omega)\geq n_1$ such that, for all $n\geq m_\vare^{*}(\omega)$, we have
 \begin{align*}
    -(C-\vare)S_n(\omega) -  T_n^{(C-\vare)}(\omega)
      &< \sqrt{\EE(D(\xi_1,t_0)^2)} - \frac{\vare}{2} \EE(-\partial_2 D(\xi_1,t_0)).
 \end{align*}
On the other hand, there exists a sequence $(n_k^{\diamond}(\omega))_{k\in\NN}$ such that
 $n_k^{\diamond}(\omega)\in\NN\cap (m_\vare^{*}(\omega),\infty)$, $n_k^{\diamond}(\omega)\uparrow \infty$ as $k\to\infty$, and
 \begin{align*}
   R_{n_k^{\diamond}(\omega)}(\omega)> \sqrt{\EE(D(\xi_1,t_0)^2)} - \frac{\vare}{2} \EE(-\partial_2 D(\xi_1,t_0)),
    \qquad k\in\NN.
 \end{align*}
Thus \eqref{help_LIL_8} holds with $n:=n_k^{\diamond}$ for all $k\in\NN$ and \eqref{help_LIL_7.5} is valid for the same values of $n$.
This implies that, for each $\vare\in(0,C)$ and $\omega\in \Omega^{(C-\vare)}$,
 \[
    \limsup_{n\to\infty} \frac{\mu_n(\omega)-t_0}{a_n} \geq  C-\vare.
 \]
By choosing $\vare:=1/m$ with $m>1/C$, $m\in\NN$,
 we have $\limsup_{n\to\infty} \frac{\mu_n(\omega)-t_0}{a_n} \geq C$ for $\omega\in \bigcap_{m=\lceil 1/C\rceil}^\infty \Omega^{(C-1/m)}$, where
 $\PP(\bigcap_{m=\lceil 1/C\rceil}^\infty \Omega^{(C-1/m)} )=1$.
This yields $\limsup_{n\to\infty} \frac{\mu_n-t_0}{a_n} \geq C$ almost surely.

All in all, we get that $\limsup_{n\to\infty} \frac{\mu_n-t_0}{a_n} = C$ almost surely, yielding \eqref{help_LIL_1}.

Finally, we turn to derive \eqref{help_LIL_2} by applying \eqref{help_LIL_1} and a ``reflection'' method.

Let $D^*:(-I)^2\to\RR$, $D^*(x,t):=-D(-x,-t)$, $x,t\in -I$.
Then it immediately follows that $D^*$ is a deviation.
For each $n\in\NN$, $x_1,\ldots,x_n\in -I$, the equation $\sum_{i=1}^n D^*(x_i,t)=0$, $t\in -I$,
 is equivalent to $\sum_{i=1}^n D(-x_i,-t)=0$, $-t\in I$,
 and hence
 \[
   \mathscr{M}_n^{D^*}(x_1,\ldots,x_n) = -\mathscr{M}_n^{D}(-x_1,\ldots,-x_n), \qquad n\in\NN, \; -x_1,\ldots,-x_n\in I.
 \]
Further, the equation $\EE(D^*(-\xi_1,t))=0$, $t\in -I$, is equivalent to $\EE(-D(\xi_1,-t))=0$, $-t\in I$, and hence we have
 $t_0^* := \mathbb{M}^{D^*}(-\xi_1) = - \mathbb{M}^D(\xi_1)=-t_0$.
It is easy to see that $-\xi_1$ satisfies the interior point condition with respect to $D^*$, and hence assumptions (i) and (ii) hold
 for $D^*$ and the random variable $-\xi_1$.
It is also not difficult to check that assumption (iii') is valid for $D^*$, for the  random variable $-\xi_1$ with $H^*(x):=H(-x)$ for $x\in -I$, $h^*(s):=h(s)$, $s\in(0,r)$, and $t_0^*=\mathbb{M}^{D^*}(-\xi_1)$.
Furthermore, we have that $\EE(D^*(-\xi_1,t_0^*)^2)=\EE(D(\xi_1,t_0)^2)$ and
 $\EE(\partial_2 D^*(-\xi_1,t_0^*))=\EE(\partial_2 D(\xi_1,t_0))$, which implies that assumption (iv) is also valid for this transformed setting.

Using \eqref{help_LIL_1} for the i.i.d.\ sequence of random variables $(-\xi_n)_{n\in\NN}$, we have
 \begin{align}\label{help_LIL_11}
   \limsup_{n\to\infty} \frac{\mathscr{M}_n^{D^*}(-\xi_1,\ldots,-\xi_n) - t_0^*}{\sqrt{2n^{-1}\ln(\ln(n))}}
       \ase  \frac{\sqrt{\EE(D^*(-\xi_1,t_0^*)^2)}}{ \EE(-\partial_2 D^*(-\xi_1,t_0^*))},
 \end{align}
Consequently, by \eqref{help_LIL_11} and the previous discussion, we obtain
 \begin{align*}
   \limsup_{n\to\infty} \frac{-\mathscr{M}_n^D(\xi_1,\ldots,\xi_n) - (-t_0)}{\sqrt{2n^{-1}\ln(\ln(n))}}
       \ase  \frac{\sqrt{\EE( D(\xi_1,t_0)^2 )}}{ \EE(-\partial_2 D(\xi_1,t_0))}
             = C,
 \end{align*}
which immediately implies \eqref{help_LIL_2}, as desired.
\proofend

\section{Large deviations for deviation means}
\label{section_LD}

We formulate a large deviation result for deviation means of i.i.d.\ random variables.

\begin{Thm}\label{Thm_Dev_mean_large_dev}
Let $I\subset\RR$ be a nondegenerate interval, let $D:I^2\to\RR$ be a deviation and $(\xi_n)_{n\in\NN}$ be a sequence of
 i.i.d.\ random variables.
Assume that
 \vspace{-12pt}
 \begin{enumerate}[(i)]\itemsep=-4pt
 \item $D$ is Borel measurable in the first variable;
 \item $\PP(\xi_1\in I)=1$ and for each $x\in I$, the random variable $D(\xi_1,x)$ \ is nondegenerate;
 \item for each $c>0$ and $x\in I$, we have
       \[
        \varphi(c,x):=\EE(\ee^{c D(\xi_1,x)})<\infty;
       \]
  \item there exists a point $t_0\in I_{\xi_1}^D$ such that $\EE(D(\xi_1,t_0))=0$, i.e., $t_0=\mathbb{M}^D(\xi_1)$.
 \end{enumerate}\vspace{-12pt}
Then for each $x>t_0$, $x\in I$ with $\esssup(D(\xi_1,x))>0$, we have
 \begin{align*}
   \lim_{n\to\infty} \frac{1}{n}\ln\Big(\PP\big(\mathscr{M}_n^{D}(\xi_1,\ldots,\xi_n) \geq x \big)\Big)
        =\sup_{m\in\NN} \frac{1}{m} \ln\Big( \PP\big(\mathscr{M}_m^{D}(\xi_1,\ldots,\xi_m) \geq x \big) \Big)
        = \inf_{c>0} \ln\big(\varphi(c,x)\big).
 \end{align*}
The equalities above can be also rewritten in the following 'exponential' form
 \begin{equation}\label{E:GCT}
   \lim_{n\to\infty}  \sqrt[n]{\PP\big(\mathscr{M}_n^{D}(\xi_1,\ldots,\xi_n) \geq x \big)}
   =\sup_{m\in\NN}  \sqrt[m]{\PP\big(\mathscr{M}_m^{D}(\xi_1,\ldots,\xi_m) \geq x \big)}
    =  \inf_{c>0}\, \varphi(c,x).
   \end{equation}
\end{Thm}

\noindent{\bf Proof.}
For each $x\in\RR$, the following equalities hold:
  \begin{align}\label{E:Aux}
    \PP\Big(\mathscr{M}_n^D(\xi_1,\dots,\xi_n) \geq  x \Big)
     = \PP\left(\sum_{i=1}^n D(\xi_i,x) \geq  0 \right)
     =\PP\left(\frac1n\sum_{i=1}^n D(\xi_i,x) \geq  0 \right).
 \end{align}
According to assumption (iii), for all $x\in I$, we have that the positive part $[D(\xi_1,x)]^+$ of $D(\xi_1,x)$ is integrable,
 i.e., $\EE([D(\xi_1,x)]^+)<\infty$, and therefore, $\EE(D(\xi_1,x))$ is well-defined and belongs to $[-\infty,\infty)$.
For each $x>t_0$, $x\in I$, by Lemma \ref{Lem_Dev_mean_random_1}, we have $\EE(D(\xi_1,x)) < \EE(D(\xi_1,t_0))=0$.
Hence for each $x>t_0$, $x\in I$ with $\esssup(D(\xi_1,x))>0$,
we can apply the exponential form of Cram\'er's theorem on large deviations for the sequence of i.i.d.\ random variables $D(\xi_i,x)$, $i\in\NN$ and $y:=0$ (see Theorem \ref{Thm_Cramer}) and we obtain
 \begin{align}\label{help3}
  \lim_{n\to\infty} \sqrt[n]{\PP\left( \frac{1}{n}\sum_{i=1}^n D(\xi_i,x) \geq 0 \right)}
  =\sup_{m\in\NN} \sqrt[m]{\PP\left( \frac{1}{m}\sum_{i=1}^m D(\xi_i,x) \geq 0 \right)}
      =\inf_{c>0} \varphi(c,x),
 \end{align}
which, in view of \eqref{E:Aux}, yields the equality \eqref{E:GCT}. This shows that the first equation of the assertion is also valid.
\proofend

\section{Examples for limit theorems}
\label{section_ELT}

We specialize Theorems \ref{Thm_SLLN_deviation}, \ref{Thm_CLT_deviation} and \ref{Thm_LIL_deviation} to Bajraktarevi\'c means by deriving a strong law of large numbers, a central limit theorem and a law of the iterated logarithm for these means.
The strong law of large numbers and the central limit theorem in question have already been proved in Barczy and Burai \cite[Theorem 2.1]{BarBur}, however, the law of the iterated logarithm in this setting is a new result, and according to our knowledge, it is new even for quasi-arithmetic means.

\begin{Ex}[Bajraktarevi\'c means]\label{Ex_3_Bajraktarevic}
Let $I$ be a nondegenerate interval of $\RR$, let $f:I\to\RR$ be a continuous and strictly increasing function and let $p:I\to (0,\infty)$ be a Borel measurable function. Let $(\xi_n)_{n\in\NN}$ be a sequence of i.i.d.\ random variables such that $\PP(\xi_1\in I)=1$, $\EE(p(\xi_1))<\infty$ and
 $\EE(p(\xi_1)\vert f(\xi_1)\vert)<\infty$.
As it was mentioned in the introduction, for each $n\in\NN$, an $n$-variable Bajraktarevi\'c mean corresponding to $f$ and $p$ is nothing else but an $n$-variable deviation mean corresponding to the deviation $D:I^2\to\RR$, $D(x,t):=p(x)(f(x)-f(t))$, $x,t\in I$.
By part (ii) of Example \ref{Ex_1}, $D$ is Borel measurable in the first variable and $I^D_\xi = I$.
Thus assumption (i) of Theorem \ref{Thm_SLLN_deviation} is satisfied.

Further, the equation $\EE(D(\xi_1,t))=0$, $t\in I$, takes the form
 \[
   0 = \EE(p(\xi_1)f(\xi_1)) - f(t)\EE(p(\xi_1)), \qquad t\in I,
 \]
and hence
 \[
    f(t) = \frac{\EE(p(\xi_1)f(\xi_1))}{\EE(p(\xi_1))},\qquad t\in I.
 \]
Using that $f$ is strictly monotone and Corollary \ref{Cor_1}, we have
 \[
    t_0:=\mathbb{M}^D(\xi_1) = f^{-1}\left( \frac{\EE(p(\xi_1)f(\xi_1))}{\EE(p(\xi_1))}\right) \in I.
 \]

Supposing that $\mathbb{M}^D(\xi_1)$ belongs to the interior of $I$, we have assumption (ii) of Theorem \ref{Thm_SLLN_deviation} is satisfied as well, and hence
  Theorem \ref{Thm_SLLN_deviation} yields
\[
 \mathscr{B}^{f,p}_n(\xi_1,\ldots,\xi_n) \as f^{-1}\left( \frac{\EE(p(\xi_1)f(\xi_1))}{\EE(p(\xi_1))} \right)  \qquad \text{as $n\to\infty$,}
\]
as we expected (compare with Theorem \ref{Thm_CLT_Baj_mean}). Here we call the attention to the fact that Barczy and Burai \cite[Theorem 2.1]{BarBur} supposed that $f(I)$ is closed for deriving a strong law of large numbers for Bajraktarevi\'c means, however, as we have seen, it is a superfluous assumption in their theorem.

In what follows, additionally, let us suppose that $\EE(p(\xi_1)^2)<\infty$, $\EE(p(\xi_1)^2 f(\xi_1)^2)<\infty$,
 $\xi_1$ is not a constant random variable (equivalently, $\PP(\xi_1=t_0)<1$) and $f$ is differentiable at $t_0$ with a nonzero derivative.
Next, we check the assumptions of Theorem \ref{Thm_CLT_deviation}.
In fact, we have already verified that the assumptions (i) and (ii) of Theorem \ref{Thm_CLT_deviation} hold (see the discussion above). By the differentiability of $f$ at $t_0$, it follows that the function $h:(0,r)\to\RR$ defined by
\[
  h(s):=\max\bigg(\bigg|\frac{f(t_0+s)-f(t_0)}{s}-f'(t_0)\bigg|,
  \bigg|\frac{f(t_0)-f(t_0-s)}{s}-f'(t_0)\bigg|\bigg), \qquad s\in(0,r),
\]
(where $r>0$ is chosen such that $(t_0-r,t_0+r)\subset I$) possesses the limit property $h(s)\to0$ as $s\to0^+$.
Therefore, the assumption (iii) of Theorem \ref{Thm_CLT_deviation} is valid with $\alpha=0$, $H(x):=p(x)$, $x\in I$, since
 \begin{align*}
   \widetilde f(t):= \left\vert \frac{f(t)-f(t_0)}{t-t_0} - f'(t_0) \right\vert
                 & \leq \max\big(\widetilde f(t), \widetilde f(2t_0-t)\big)
                   = \max\big(\widetilde f(t_0+\vert t-t_0\vert), \widetilde f(t_0-\vert t-t_0\vert)\big)\\
                 & = h(\vert t-t_0\vert),\qquad t\in(t_0-r,t_0+r), \; t\ne t_0.
 \end{align*}
In this case $q=1$, and $\EE(H(\xi_1)^q)<\infty$ follows from $\EE(p(\xi_1))<\infty$.
Finally,
 \begin{align*}
    \EE(D(\xi_1,t_0)^2)
      = \EE(p(\xi_1)^2(f(\xi_1) - f(t_0))^2)
      \leq 2\EE(p(\xi_1)^2 f(\xi_1)^2) + 2 f(t_0)^2 \EE(p(\xi_1)^2)<\infty,
 \end{align*}
and $\EE(D(\xi_1,t_0)^2)=0$ holds if and only if $\PP(f(\xi_1)= f(t_0))=1$, and, equivalently, using the strict monotonicity of $f$, if and only if $\PP(\xi_1= t_0)=1$, and, since we have assumed that $\xi_1$ is not a constant random variable, we can conclude that $\EE(D(\xi_1,t_0)^2)>0$.
Moreover, $\EE(-\partial_2D(\xi_1,t_0)) = \EE(p(\xi_1)f'(t_0))=f'(t_0)\EE(p(\xi_1))\in(0,\infty)$, and hence the assumption (iv) of Theorem \ref{Thm_CLT_deviation} holds, too.

Thus one can apply Theorem \ref{Thm_CLT_deviation} yielding
  \[
  \sqrt{n}( \mathscr{B}^{f,p}_n(\xi_1,\ldots,\xi_n) - t_0)\distr \cN\left(0,\frac{\EE( D(\xi_1,t_0)^2 )}{( \EE(\partial_2 D(\xi_1,t_0)) )^2}\right)
             \qquad \text{as $n\to\infty$.}
  \]
Here
 \begin{align*}
   \EE( D(\xi_1,t_0)^2 )
     &= \EE( (p(\xi_1)(f(\xi_1) - f(t_0)) )^2 )
      = \EE\left(  (p(\xi_1))^2 \left( f(\xi_1) - \frac{\EE(p(\xi_1)f(\xi_1))}{\EE(p(\xi_1))} \right)^2 \right)\\
     &= \frac{1}{(\EE( p(\xi_1)))^2}
        \Big[ (\EE( p(\xi_1) ))^2 \var(p(\xi_1)f(\xi_1))
        + ( \EE( p(\xi_1) f(\xi_1) ) )^2 \var(p(\xi_1))\\
    &\phantom{= \frac{1}{(\EE( p(\xi_1)))^2} \Big[}
       - 2 \EE(p(\xi_1)) \EE(p(\xi_1)f(\xi_1))\cov(p(\xi_1),p(\xi_1)f(\xi_1))\Big] ,
 \end{align*}
 and
 \begin{align*}
   \big( \EE( \partial_2 D(\xi_1,t_0) ) \big)^2
    = \big( \EE( -p(\xi_1) f'(t_0) ) \big)^2
    = \left( \EE( p(\xi_1)) \right)^2 \left(f'\left( f^{-1}\left( \frac{\EE(p(\xi_1)f(\xi_1))}{\EE(p(\xi_1))} \right) \right)\right)^2 .
 \end{align*}
This yields that
 \[
 \frac{\EE( D(\xi_1,t_0)^2 )}{( \EE(\partial_2 D(\xi_1,t_0)))^2}
  = \sigma_{f,p}^2,
 \]
 where $\sigma_{f,p}^2$ is given in Theorem \ref{Thm_CLT_Baj_mean}, and hence Theorem \ref{Thm_CLT_deviation} yields \eqref{help7_Baj_mean}, as we expected
 (compare with Theorem \ref{Thm_CLT_Baj_mean}).

Now we turn to establish the law of the iterated logarithm for Bajraktarevi\'c means of i.i.d.\ random variables.
Note that in case of $\alpha=0$, the assumption (iii') of Theorem \ref{Thm_LIL_deviation} coincides with the assumption (iii) of Theorem \ref{Thm_CLT_deviation}.
Consequently, all the assumptions of Theorem \ref{Thm_LIL_deviation} hold and we have
 \begin{align*}
    \limsup_{n\to\infty} \frac{\mathscr{B}^{f,p}_n(\xi_1,\ldots,\xi_n) - t_0}{\sqrt{2n^{-1}\ln(\ln(n))}}
       \ase  \sqrt{\sigma_{f,p}^2}
 \end{align*}
 and
 \begin{align*}
    \liminf_{n\to\infty} \frac{\mathscr{B}^{f,p}_n(\xi_1,\ldots,\xi_n) - t_0}{\sqrt{2n^{-1}\ln(\ln(n))}}
       \ase  -\sqrt{\sigma_{f,p}^2}.
 \end{align*}
In the special case $p(x)=1$, $x\in I$, i.e., in case of a quasi-arithmetic mean with a generator $f$, we have
 \begin{align*}
    \limsup_{n\to\infty} \frac{\mathscr{A}^f_n(\xi_1,\ldots,\xi_n) - \mathbb{A}^f(\xi_1) }{\sqrt{2n^{-1}\ln(\ln(n))}}
           \ase  \frac{\sqrt{\var(f(\xi_1))}}{f'(\mathbb{A}^f(\xi_1))},
 \end{align*}
 and
 \begin{align*}
     \liminf_{n\to\infty} \frac{\mathscr{A}^f_n(\xi_1,\ldots,\xi_n) - \mathbb{A}^f(\xi_1) }{\sqrt{2n^{-1}\ln(\ln(n))}}
           \ase  -\frac{\sqrt{\var(f(\xi_1))}}{f'(\mathbb{A}^f(\xi_1))},
 \end{align*}
 where we recall that $\mathscr{A}^f_n(\xi_1,\ldots,\xi_n)$ denotes the quasi-arithmetic mean of $\xi_1,\ldots,\xi_n$,
 and $\mathbb{A}^f(\xi_1) = f^{-1}(\EE(f(\xi_1)))$ is the quasi-arithmetic expected value of $\xi_1$.
According to our knowledge, laws of iterated logarithms have not been established for Bajraktarevi\'c and quasi-arithmetic means so far,
 and hence these results seem to be new as well.
\proofend
\end{Ex}

In what follows, we provide an example for a deviation mean which is not a Bajraktarevi\'c mean, and then we specialize Theorems \ref{Thm_SLLN_deviation}, \ref{Thm_CLT_deviation} and \ref{Thm_LIL_deviation} to this particular deviation mean.

\begin{claim}\label{Cla_1}
Let $I:=(0,\infty)$ and $D:I^2\to\RR$, $D(x,t):=x(x-t)+x^2-t^2$, $x,t\in I$.
Then $D$ is a deviation which is Borel measurable in the first variable,
 and the family of deviation means $\{\mathscr{M}^D_n, n\in\NN\}$ cannot be a family of Bajraktarevi\'c means, i.e., there do not exist a continuous and strictly monotone function $f:I\to\RR$ and a function $p: I\to (0,\infty)$ such that $\{\mathscr{B}^{f,p}_n, n\in\NN\}$ coincides with $\{\mathscr{M}^D_n, n\in\NN\}$.
\end{claim}

\noindent{\bf Proof.}
Since $\partial_2 D(x,t) = -x -2t<0$, $x,t\in I$, and $D(t,t)=0$, $t\in I$, we have that $D$ is a deviation.
On the contrary, let us suppose that $\{\mathscr{M}^D_n, n\in\NN\}$ is a family of Bajraktarevi\'c means.
Then, using a result of Dar\'oczy and P\'ales \cite{DarPal82} (for more details, see the end of Remark \ref{Rem1}),
 $D$ can be written in the form $D(x,t) = d(t)p(x)(f(x)-f(t))$, $x,t\in I$, where $f:I\to\RR$ is a continuous and strictly increasing function,
 $p:I\to\RR$ is a positive function and $d:I\to\RR$ is a positive continuous function.
Consequently, $D$ satisfies the conditions and part (i) of the Characterization Theorem 1 of P\'ales \cite{Pal87}, which says that, for all $n\in\NN$, $x_1,\dots,x_n\in I$ and $\lambda_1,\dots,\lambda_n\in\RR$, the set
\[
  \Bigg\{ t\in I :  \sum_{i=1}^n \lambda_iD(x_i,t)\leq 0 \Bigg\}
\]
 must be convex (i.e., it must the empty set, a singleton or the interval $I$).
Here
 \begin{align*}
   \Bigg\{ t\in I :  \sum_{i=1}^n \lambda_iD(x_i,t)\leq 0 \Bigg\}
   & = \Bigg\{ t\in I : \sum_{i=1}^n \lambda_i (2x_i^2-x_it-t^2)\leq 0\Bigg\} \\
   & = \Bigg\{ t\in I :  -\left(\sum_{i=1}^n \lambda_i \right)t^2 -\left(\sum_{i=1}^n\lambda_i x_i \right)t
           + 2 \sum_{i=1}^n \lambda_i x_i^2 \leq 0 \Bigg\}.
 \end{align*}
Hence, in order to get a contradiction, it is enough to find $n\in\NN$, $x_1,\dots,x_n\in I$ and $\lambda_1,\dots,\lambda_n\in\RR$ with $\sum_{i=1}^n \lambda_i>0$ such that the equation $-\left(\sum_{i=1}^n \lambda_i \right)t^2 -\left(\sum_{i=1}^n\lambda_i x_i \right)t + 2 \sum_{i=1}^n \lambda_i x_i^2 = 0$ has two positive distinct solutions, which take the forms
 \[
    t_{\pm}:=\frac{-\sum_{i=1}^n\lambda_ix_i\pm\sqrt{(\sum_{i=1}^n\lambda_ix_i)^2+8\sum_{i=1}^n\lambda_i\sum_{i=1}^n\lambda_ix_i^2}}{2\sum_{i=1}^n\lambda_i}.
 \]
To give such examples, let $n:=2$, $\lambda_1:=51$, $\lambda_2:=-50$, $x_1:=50$, and $x_2:=101$.
Then
 \begin{align*}
 &\sum_{i=1}^2\lambda_i=1,\qquad \sum_{i=1}^2\lambda_ix_i=-2500, \qquad
  \sum_{i=1}^2\lambda_ix_i^2=-382550,\\
 &\left(\sum_{i=1}^2\lambda_ix_i\right)^2+8\sum_{i=1}^2\lambda_i\sum_{i=1}^2\lambda_ix_i^2 = 3189600,
 \end{align*}
 and hence
 \[
 t_\pm=\frac{2500\pm\sqrt{3189600}}{2}.
 \]
Consequently,
 \[
    \left\{ t\in I :  \sum_{i=1}^2 \lambda_iD(x_i,t)\leq 0 \right\}
      = (0,t_-]\cup [t_+,\infty),
 \]
 where $t_-\approx 357,03$ and $t_+\approx 2142,97$, which is not a convex set, yielding a contradiction.
\proofend

\begin{Ex}\label{Ex_4_nonBajraktarevic}
Let $I:=(0,\infty)$ and $D:I^2\to\RR$, $D(x,t):=x(x-t)+x^2-t^2$, $x,t\in I$.
Then, by Claim \ref{Cla_1}, $D$ is a deviation which is Borel measurable in the first variable, and the family of deviation means $\{\mathscr{M}^D_n, n\in\NN\}$ cannot be a family of Bajraktarevi\'c means.
Further, for each $n\in\NN$ and $x_1,\ldots,x_n\in I$, the equation $\sum_{i=1}^n D(x_i,t) = 0$, $t\in I$, takes the form
 \[
    2 \sum_{i=1}^n x_i^2 - \left( \sum_{i=1}^n x_i \right)t-nt^2 = 0,\qquad t>0,
 \]
 and it has a unique solution
 \begin{align*}
   \mathscr{M}_n^{D}(x_1,\ldots,x_n)
      & = -\frac{1}{2} \frac{\sum_{i=1}^n x_i}{n} + \sqrt{\frac{1}{4}\left(\frac{\sum_{i=1}^n x_i}{n}\right)^2  + 2 \frac{\sum_{i=1}^n x_i^2}{n}}\\
      & = -\frac{1}{2} \mathscr{A}^f_n(x_1,\ldots,x_n) + \sqrt{\frac14\big(\mathscr{A}^f_n(x_1,\ldots,x_n) \big)^2  + 2 \big(\mathscr{A}^g_n(x_1,\ldots,x_n)\big)^2 },
 \end{align*}
 where $f(x):=x$, $x>0$, $g(x):=x^2$, $x>0$, and we used the notations for quasi-arithmetic means introduced in Definition \ref{Def_quasi_arithmetic}.
Due to the well-known ineqaulity $\mathscr{A}^f_n(x_1,\ldots,x_n)\leq \mathscr{A}^g_n(x_1,\ldots,x_n)$,
  it follows that $\mathscr{M}_n^{D}(x_1,\ldots,x_n)\geq \mathscr{A}^f_n(x_1,\ldots,x_n)$ for each $n\in\NN$ and $x_1,.\ldots,x_n\in I$.

Let $(\xi_n)_{n\in\NN}$ be a sequence of i.i.d.\ random variables such that $\PP(\xi_1\in I)=1$ and $\EE(\xi_1^2)<\infty$.
Then for each $t\in I$, we have
 \begin{align*}
    \EE(\vert D(\xi_1,t)\vert) = \EE(\vert 2\xi_1^2 - \xi_1 t - t^2\vert)
                               \leq 2\EE(\xi_1^2) + t\EE(\xi_1) + t^2<\infty,
 \end{align*}
 and hence $I^D_{\xi_1} = I$.
Further, the equation $\EE(D(\xi_1,t))=0$, $t\in I$, takes the form
 \[
  \EE(2\xi_1^2 - \xi_1 t - t^2)=0,\qquad t>0,
 \]
 which is equivalent to $t^2 + \EE(\xi_1)t - 2\EE(\xi_1^2)=0$, $t>0$.
Consequently,
 \[
    t_0:=\mathbb{M}^D(\xi_1) = \frac{1}{2}\Big(-\EE(\xi_1)+\sqrt{8\EE(\xi_1^2) + (\EE(\xi_1))^2}\Big)
    = -\frac{1}{2}\EE(\xi_1) + \sqrt{\frac14(\EE(\xi_1))^2+2 \EE(\xi_1^2)} .
 \]
Indeed, by Lemma \ref{Lem_Exp_conv_hull}, $\EE(\xi_1^2)\in I=(0,\infty)$,
 and hence $-\EE(\xi_1)+\sqrt{8\EE(\xi_1^2) + (\EE(\xi_1))^2} > -\EE(\xi_1)+\EE(\xi_1)=0$.
Consequently, $\mathbb{M}^D(\xi_1)\in I$, and thus $\xi_1$ satisfies the interior point condition with respect to $D$.
All in all, the assumptions (i) and (ii) of Theorem \ref{Thm_SLLN_deviation} are satisfied, and hence we obtain
 \begin{align*}
   \mathscr{M}_n^{D}(\xi_1,\ldots,\xi_n)
      &= -\frac{1}{2} \mathscr{A}^f_n(x_1,\ldots,x_n) + \sqrt{\frac{1}{4} \big(\mathscr{A}^f_n(x_1,\ldots,x_n) \big)^2  + 2 (\mathscr{A}^g_n(x_1,\ldots,x_n))^2 }\\
      & \as -\frac{1}{2}\EE(\xi_1) + \sqrt{  \frac14(\EE(\xi_1))^2 + 2 \EE(\xi_1^2) }
 \end{align*}
 as $n\to\infty$.
Of course, to derive the almost sure convergence above, we do not need to use Theorem \ref{Thm_SLLN_deviation}, it readily follows by
 Kolmogorov's strong law of large numbers for  the i.i.d.\ random variables $(\xi_n)_{n\in\NN}$ and $(\xi_n^2)_{n\in\NN}$, respectively.

In what follows, additionally, let us suppose that $\EE(\xi_1^4)<\infty$ and $\xi_1$ is not a constant random variable (equivalently, $\PP(\xi_1=t_0)<1$).
Next, we check the assumptions of Theorem \ref{Thm_CLT_deviation}.
 In fact, we have already verified that the assumptions (i) and (ii) of Theorem  \ref{Thm_CLT_deviation} hold (see the discussion above).
Further, for each $x,t\in I$, we have
 \[
   \vert \partial_2D(x,t) - \partial_2D(x,t_0) \vert
      =  \vert -x-2t - (-x-2t_0) \vert
      = 2  \vert t-t_0 \vert,
 \]
and hence the inequality \eqref{help_assump_iv_ver1} holds with $H(x)=1$, $x\in I$, and $h(s) = 2s$, $s\in[0,r)$, where $r>0$ is such that $t_0>r$ (such an $r>0$ exists, since $t_0>0$).
By part (i) of Remark \ref{Rem_CLT_a(iii)_elegseges}, the assumption (iii) of Theorem \ref{Thm_CLT_deviation} holds with the given functions $h$, $H$ and $\alpha:=0$ yielding $q=1$.
Moreover,
 \[
   \EE(D(\xi_1,t_0)^2) = \EE( (2\xi_1^2 - \xi_1t_0 - t_0^2)^2 )\leq 2\big( 4\EE(\xi_1^4) + \EE(\xi_1^2)t_0^2 + t_0^4\big)<\infty,
 \]
 and $\EE(D(\xi_1,t_0)^2)=0$ is equivalent to $\PP(2\xi_1^2 - \xi_1t_0 - t_0^2=0)=1$, and, since $2\xi_1^2 - \xi_1t_0 - t_0^2 = (\xi_1-t_0)(2\xi_1+t_0)$,
 it is equivalent to $\PP(\xi_1=t_0)=1$, \ which we excluded.
Hence $0<\EE(D(\xi_1,t_0)^2)<\infty$.
Further, $\EE(-\partial_2D(\xi_1,t_0)) =\EE(\xi_1 + 2t_0) = \EE(\xi_1) + 2t_0<\infty$ and $\EE(-\partial_2D(\xi_1,t_0))=0$ holds if and only if
 $\EE(\xi_1)=-2t_0$, which cannot happen, since $\EE(\xi_1)>0$ and $-2t_0<0$.
Hence $0<\EE(-\partial_2D(\xi_1,t_0))<\infty$, and then the assumption (iv) of Theorem \ref{Thm_CLT_deviation} is satisfied.
Thus one can apply Theorem \ref{Thm_CLT_deviation} yielding
 \[
  \sqrt{n}( \mathscr{M}_n^{D}(\xi_1,\ldots,\xi_n)  - t_0)\distr \cN\left(0,\frac{\EE( D(\xi_1,t_0)^2 )}{( \EE(\partial_2 D(\xi_1,t_0)) )^2}\right)
             \qquad \text{as $n\to\infty$.}
 \]
Here
 \[
   ( \EE(\partial_2 D(\xi_1,t_0)) )^2 = (\EE(\xi_1) + 2t_0)^2
                                      = 8\EE(\xi_1^2) + (\EE(\xi_1))^2 ,
 \]
 and, using the equality $t_0^2 + \EE(\xi_1)t_0 - 2\EE(\xi_1^2)=0$ in each step, we obtain
 \begin{align*}
  \EE( D(\xi_1,t_0)^2 )
   &= \EE( (2\xi_1^2 - \xi_1t_0 -t_0^2)^2 )
    = t_0^4 +2\EE(\xi_1)t_0^3 - 3 \EE(\xi_1^2)t_0^2 -4\EE(\xi_1^3)t_0 + 4\EE(\xi_1^4)\\
   &= \EE(\xi_1) t_0^3  - \EE(\xi_1^2) t_0^2 - 4\EE(\xi_1^3) t_0  + 4\EE(\xi_1^4)\\
   &= - (\EE(\xi_1^2)+(\EE(\xi_1))^2)t_0^2+ (2\EE(\xi_1^2)\EE(\xi_1) - 4\EE(\xi_1^3)) t_0  + 4\EE(\xi_1^4)\\
   &= ((\EE(\xi_1))^3+3\EE(\xi_1^2)\EE(\xi_1) - 4\EE(\xi_1^3)) t_0  + 4\EE(\xi_1^4)- 2(\EE(\xi_1^2))^2-2\EE(\xi_1^2)(\EE(\xi_1))^2.
 \end{align*}

Now we turn to establish the law of the iterated logarithm for the given deviation means of i.i.d.\ random variables.
Note that in case of $\alpha=0$, the assumption (iii') of Theorem \ref{Thm_LIL_deviation} coincides with the assumption (iii) of Theorem \ref{Thm_CLT_deviation}.
Consequently, Theorem \ref{Thm_LIL_deviation} yields that \eqref{help_LIL_1} and \eqref{help_LIL_2} hold with $t_0$ given above.
\proofend
\end{Ex}

The investigation of limit theorems for means of i.i.d.\ random variables possibly could be extended to further classes of means that are not deviation means, see the next remark.

\begin{Rem}
We note that Barczy and Burai \cite{BarBur} derived strong laws of large numbers and central limit theorems for so-called Cauchy quotient means, for example for the exponential Cauchy quotient mean (also called Beta-type mean) defined by
 \[
   \cB_n(x_1,\ldots,x_n):=\sqrt[n-1]{\frac{nx_1\cdots x_n}{x_1+\cdots+x_n}},
 \]
 where $n\geq 2$, $n\in\NN$ and $x_1,\ldots,x_n>0$.
However, Theorems \ref{Thm_SLLN_deviation} and Theorem \ref{Thm_CLT_deviation} cannot be applied to exponential Cauchy quotient means,
 since there does not exist a deviation $D:(0,\infty)^2\to\RR$ such that the family $\{\cB_n, n\geq2, n\in\NN \}$ coincides with the family
 $\{\mathscr{M}_n^{D}: n\geq 2, n\in\NN\}$.
Indeed, if $D:(0,\infty)^2\to\RR$ were a deviation such that the two families in question coincide, then, by the so-called repetition invariance property of
 deviation means, we would have $\cB_{nk}(x_1,\ldots,x_n,\ldots,x_1,\ldots,x_n) = \cB_n(x_1,\ldots,x_n)$ for each $n\geq 2$, $n,k\in\NN$ and $x_1,\ldots,x_n>0$,
 that is,
 \[
  \sqrt[nk-1]{\frac{nk (x_1\cdots x_n)^k}{k(x_1+\cdots+x_n)}}
    = \sqrt[n-1]{\frac{nx_1\cdots x_n}{x_1+\cdots+x_n}}, \qquad n\geq 2, \; n,k\in\NN, \;\; x_1,\ldots,x_n>0,
 \]
 which is equivalent to
 \[
  \left(\frac{\cG_n(x_1,\ldots,x_n)^{nk}}{\cA_n(x_1,\ldots,x_n)}\right)^{n-1}
      = \left(\frac{\cG_n(x_1,\ldots,x_n)^n}{\cA_n(x_1,\ldots,x_n)}\right)^{nk-1},
      \qquad n\geq 2, \; n,k\in\NN,\;\;  x_1,\ldots,x_n>0,
 \]
 where $\cA_n(x_1,\ldots,x_n):=n^{-1}(x_1+\cdots+x_n)$ and $\cG_n(x_1,\ldots,x_n):=\sqrt[n]{x_1\cdots x_n}$
 are the arithmetic and geometric means of $x_1,\ldots,x_n>0$, respectively.
By an easy calculation, it yields that $\cA_n(x_1,\ldots,x_n) = \cG_n(x_1,\ldots,x_n)$ should hold for each $n\in\NN$ and $x_1,\ldots, x_n>0$,
 which leads us to a contradiction.
Finally, we mention that a law of the iterated logarithm is not known for exponential Cauchy means if i.i.d.\ random variables.
\proofend
\end{Rem}

\appendix

\section{Appendix: Auxiliary results}

In this part, we recall or prove some auxiliary results that are used in the proofs.

\begin{Lem}\label{Lem_Exp_conv_hull}
Let $I$ be a nondegenerate interval of $\RR$, let $f:I\to\RR$ be a continuous and strictly monotone function,
 and let $\xi$ be a random variable such that $\PP(\xi\in I)=1$ and $\EE(\vert f(\xi)\vert)<\infty$.
Then $\EE (f(\xi))\in f(I)$.
\end{Lem}

\noindent{\bf Proof.}
First note that $f(I)$ is an interval and so it is a convex set.
If $\xi_n$, $n\in\NN$, are i.i.d.\ random variables such that $\xi_1\distre\xi$, then, by Kolmogorov's strong law of large numbers,
 $\frac{1}{n}\sum_{i=1}^n f(\xi_i)\as \EE(f(\xi))$ as $n\to\infty$.
So $\EE(f(\xi))$ belongs to the closure $\overline{f(I)}$ of the interval $f(I)$.
If $f(I)$ is closed, then the statement is obvious.
In what follows, let us suppose that $f(I)$ is not closed.
Then we check that $\EE(f(\xi))\in \overline{f(I)}\setminus f(I)$ cannot hold. On the contrary, let us assume it.
Then $f(I)\neq \RR$, and $[\EE(f(\xi)),\infty)\cap f(I)=\emptyset$ or $(-\infty,\EE(f(\xi))]\cap f(I)=\emptyset$
 yielding that $\PP(f(\xi) - \EE(f(\xi)) \leq 0)=1$ or $\PP(f(\xi) - \EE(f(\xi)) \geq 0)=1$.
Since $\EE(f(\xi) - \EE(f(\xi))) = 0$, it is known that it implies $\PP(f(\xi) =  \EE(f(\xi)) )=1$
 leading us to a contradiction, since $\PP(f(\xi)\in f(I)) = 1$ and $\EE(f(\xi))\notin f(I)$.
\proofend

\begin{Cor}\label{Cor_1}
Let $I$ be a nondegenerate interval of $\RR$, let $f:I\to\RR$ be a continuous and strictly monotone
 function, and let $p:I\to (0,\infty)$ be a Borel measurable function.
Let $\xi$ be a random variable with $\PP(\xi\in I)=1$, $\EE(p(\xi))<\infty$ and $\EE(p(\xi)|f(\xi)|)<\infty$.
Then
  \[
    \frac{\EE(p(\xi)f(\xi))}{\EE(p(\xi))}\in f(I).
 \]
\end{Cor}

\noindent{\bf Proof.}
Let us define the probability measure $\QQ$ on the measurable space $(\Omega,\cF)$ by
 \[
  \QQ(A):=\int_A \frac{p(\xi)}{\EE(p(\xi))}\,\dd \PP = \frac{\EE(p(\xi)\bone_{A})}{\EE(p(\xi))},\qquad A\in\cF,
 \]
where $(\Omega,\cF,\PP)$  denotes the underlying probability space.
Then $f(\xi)$ is integrable with respect to $\QQ$, and, by denoting the expectation with respect to $\QQ$ by $\EE_\QQ$, we have
 \[
    \EE_\QQ(f(\xi)) = \frac{\EE(p(\xi)f(\xi))}{\EE(p(\xi))}.
 \]
By Lemma \ref{Lem_Exp_conv_hull}, $\EE_\QQ(f(\xi))\in f(I)$, yielding the statement.
\proofend

For the proof of the central limit theorem, we formulate a direct consequence of the strong law of large numbers due to Marcinkiewicz and Zygmund, see, e.g., Kallenberg \cite[Theorem 5.23]{Kal21}. We note that the case $q=1$ is the strong law of large numbers by Kolmogorov.

\begin{Lem}\label{Lem_MZ_SLLN}
Let $(\xi_n)_{n\in\NN}$ be a sequence of i.i.d.\ random variables and let $q\in(0,1]$. Assume that $\EE(\vert \xi_1\vert^q)<\infty$.
Then $n^{-\frac{1}{q}} \sum_{i=1}^n \xi_i$ converges as $n\to\infty$ almost surely to $\EE(\xi_1)$ if $q=1$, and to $0$ if $q\in(0,1)$.
\end{Lem}

\begin{Lem}\label{Lem_aux_dev_as}
Let $I\subset\RR$ be a nondegenerate interval, let $t_0$ be an interior point of $I$, let $D:I^2\to\RR$ be a deviation and $(\xi_n)_{n\in\NN}$ be a sequence of i.i.d.\ random variables such that $\PP(\xi_1\in I)=1$.
Assume that
 \begin{itemize}
 \item[(a)] $D$ is Borel measurable in the first variable;
 \item[(b)] $D$ is differentiable in the second variable at $t_0$, and there exist $r>0$, $\alpha\in[0,\infty)$,
           a nonnegative Borel measurable function $H:I\to\RR$ and a nonnegative function $h:(0,r)\to\RR$ such that
           \begin{align*}
           \bigg|\frac{D(x, t)- D(x, t_0)}{t-t_0}- \partial_2D(x, t_0)\bigg| \leq H(x)h(|t-t_0|)
           \end{align*}
       for $x\in I$, $t\ne t_0$, $t\in (t_0 - r, t_0 +r)\subset I$,
       and
     \[
        \limsup_{s\to0^+}\frac{h(s)}{s^\alpha}\,\,
         \begin{cases}
             =0 & \mbox{if }\, \alpha=0,\\[2mm]
             <\infty & \mbox{if }\, \alpha>0,
         \end{cases}
         \qquad\mbox{and}\qquad
          \EE(H(\xi_1)^q)<\infty,
     \]
    where $q\in(0,1]$ if $\alpha=0$, and $q\in(0,1)$ if $\alpha>0$.
\end{itemize}
Further, let $(\kappa_n)_{n\in\NN}$ and $(\tau_n)_{n\in\NN}$ be sequences of real numbers such that
 \begin{itemize}
    \item[(c)] $\kappa_n\ne0$, $n\in\NN$, $\lim_{n\to\infty}\kappa_n = 0$, and
              \[
                \limsup_{n\to\infty} n^{\frac{1}{q}} \vert \tau_n\vert \vert \kappa_n\vert^{\alpha+1} <\infty.
              \]
 \end{itemize}
Then
 \[
  \tau_n \sum_{i=1}^n \Big(D\big(\xi_i, t_0 + \kappa_n \big)- D(\xi_i, t_0)
                                        - \kappa_n \partial_2D(\xi_i, t_0)\Big)
       \as 0 \qquad \text{as \ $n\to\infty$.}
 \]
\end{Lem}

\noindent{\bf Proof.}
Since $\lim_{n\to\infty}\kappa_n = 0$, we may assume that $|\kappa_n|<r$ for $n\in\NN$.
Note that the assumptions (a) and (b) yield that $D\big(\xi_i, t_0\big)$, $D\big(\xi_i, t_0 + \kappa_n \big)$ and
 $\partial_2D(\xi_i, t_0)$, $i\in\NN$, $n\in\NN$, are random variables.
Indeed, for $D(\xi_i, t_0)$ and $D(\xi_i, t_0+\kappa_n)$, this follows from assumption (a).
Further, for each $x\in I$,
 \[
  \partial_2D(x, t_0) = \lim_{n\to\infty} \frac{D(x,t_0+\kappa_n) - D(x,t_0)}{\kappa_n},
 \]
 where, by assumption (a), the function $I\ni x\mapsto \kappa_n^{-1}(D(x,t_0+\kappa_n) - D(x,t_0))$ is Borel measurable
 for each $n\in\NN$.
Using that pointwise limit of Borel measurable functions is Borel measurable, we have the Borel measurability
 of $I\ni x\mapsto \partial_2D(x, t_0)$, yielding that $\partial_2D(\xi_i, t_0)$, $i\in\NN$ are random variables.

Using assumption (b), for each $n\in\NN$, we have that
\begin{align*}
  \bigg\vert \tau_n \sum_{i=1}^n &\Big(D\big(\xi_i, t_0 + \kappa_n \big)- D(\xi_i, t_0)
                         -\kappa_n \partial_2D(\xi_i, t_0)\bigg) \bigg\vert\leq
\end{align*}
\begin{align*}
  &\leq \vert \tau_n\kappa_n\vert \sum_{i=1}^n \left\vert \frac{D\big(\xi_i, t_0 + \kappa_n \big)- D(\xi_i, t_0)}{\kappa_n}
                                                                - \partial_2D(\xi_i, t_0) \right\vert
  \leq \vert \tau_n\kappa_n\vert \sum_{i=1}^n H(\xi_i) h(|\kappa_n|).
\end{align*}
Here
 \[
  \vert \tau_n\kappa_n\vert \sum_{i=1}^n H(\xi_i) h(|\kappa_n|)
     = n^{\frac{1}{q}} \vert\tau_n\kappa_n\vert \vert\kappa_n\vert^\alpha \cdot \frac{h(\vert\kappa_n\vert)}{\vert\kappa_n\vert^\alpha}
       \cdot \frac{1}{n^{\frac{1}{q}}} \sum_{i=1}^n H(\xi_i),
 \]
 where $q$ was given in assumption (b).
Therefore, using that $h$ and $H$ are nonnegative, we have
 \begin{align}\label{help_lemma_as_1}
  \begin{split}
    & \limsup_{n\to\infty}
      \vert \tau_n\kappa_n\vert \sum_{i=1}^n H(\xi_i) h(|\kappa_n|) \\
    &\qquad \leq \limsup_{n\to\infty} n^{\frac{1}{q}} \vert\tau_n\vert \vert\kappa_n\vert^{\alpha+1}\cdot
            \limsup_{n\to\infty} \frac{h(\vert\kappa_n\vert)}{\vert\kappa_n\vert^\alpha} \cdot
            \limsup_{n\to\infty} \frac{1}{n^{\frac{1}{q}}} \sum_{i=1}^n H(\xi_i).
  \end{split}
 \end{align}

If $\alpha=0$, then $q\in(0,1]$. In this case, by assumption (b), the second $\limsup$ on the right hand side of \eqref{help_lemma_as_1} equals zero and,
 by Lemma \ref{Lem_MZ_SLLN}, we have that $n^{-\frac{1}{q}} \sum_{i=1}^n H(\xi_i)$ converges to $\EE(H(\xi_1))$ (in case $q=1$) or to $0$
 (in case $q\in(0,1)$) as $n\to\infty$ almost surely.
Further, by assumption (c), the first $\limsup$ on the right hand side of \eqref{help_lemma_as_1} is finite.

If $\alpha>0$, then $q\in(0,1)$. In this case, by assumption (b), the second $\limsup$ on the right hand side of \eqref{help_lemma_as_1} is finite
 and, by Lemma \ref{Lem_MZ_SLLN}, we have that $n^{-\frac{1}{q}} \sum_{i=1}^n H(\xi_i)$ converges to $0$ as $n\to\infty$ almost surely.
Further, by assumption (c), the first $\limsup$ on the right hand side of \eqref{help_lemma_as_1} is finite.

Consequently, in both cases (i.e., for $\alpha\in[0,\infty)$), we have
\[
   \vert \tau_n\kappa_n\vert \sum_{i=1}^n H(\xi_i) h(|\kappa_n|) \as 0 \qquad \hbox{as }\, n\to\infty,
 \]
 which yields the statement.
\proofend

Next, we recall Cram\'er's theorem on large deviations, cf.\ Athreya and Lahiri \cite[Theorem 11.4.6]{AthLah}
 and Cerf and Petit \cite{CerPet}.

\begin{Thm}[Cram\'er's theorem]\label{Thm_Cramer}
Let $(\xi_n)_{n\in\NN}$ be a sequence of i.i.d.\ nondegenerate random variables such that
 \[
    \phi(c):=\EE(\ee^{c\xi_1})<\infty, \qquad c\in\RR_{++}.
 \]
Then for each $y\in\big(\EE(\xi_1),\esssup(\xi_1)\big)$, we have
 \begin{align*}
   \lim_{n\to\infty} \frac{1}{n}\ln\Big( \PP\Big( \frac{1}{n}\sum_{i=1}^n \xi_i \geq y \Big)\Big)
       = \sup_{m\in\NN} \frac{1}{m}\ln\Big( \PP\Big( \frac{1}{m}\sum_{i=1}^m \xi_i \geq y \Big)\Big)
       = -\gamma(y),
 \end{align*}
 where $\gamma(y):=\sup_{c>0}\{cy - \ln(\phi(c))\}$, $y\in\RR$,
 is the Fenchel--Legendre transform of $\ln\circ\phi$.
The equalities above can be also rewritten in the  following 'exponential'  form
 \begin{align*}
   \lim_{n\to\infty} \sqrt[n]{ \PP\Big( \frac{1}{n}\sum\nolimits_{i=1}^n \xi_i \geq y \Big)}
   =\sup_{m\in\NN} \sqrt[m]{ \PP\Big( \frac{1}{m}\sum\nolimits_{i=1}^m \xi_i \geq y \Big)}
       = \inf_{c>0} \phi(c)\ee^{-cy}.
 \end{align*}
\end{Thm}

\section*{Data availability statements}
Data sharing is not applicable to this article as no datasets were generated or analyzed during the current study.

\section*{Acknowledgements}
We would like to thank the referee for the comments that helped us improve the paper.

\bibliographystyle{plain}
\bibliography{deviation_means_limit_thm_bib}

\end{document}